\documentclass[aip,amsmath,amssymb,reprint]{revtex4-1}

\usepackage{graphicx}% Include figure files
\usepackage{dcolumn}% Align table columns on decimal point
\usepackage{bm}% bold math
\usepackage[utf8]{inputenc}
\usepackage[T1]{fontenc}
\usepackage{mathptmx}
\usepackage{verbatim}
\usepackage{fdsymbol}
\usepackage{graphicx}
\usepackage{sidecap}
\usepackage{placeins}
\usepackage{xcolor}
\usepackage{amsfonts,amsthm,amsmath,amssymb,amscd,graphicx}
\usepackage{mathtools}
\usepackage{hyperref}
\usepackage{amsmath}

\newtheorem{theorem}{Theorem}

\makeatletter
\DeclareTextCommand{\textprime}{\encodingdefault}{%
  \mbox{$\m@th'\kern-\scriptspace$}% 
}
\makeatother

\begin{document}

\preprint{AIP/123-QED}

\title{Ordered intricacy of Shilnikov saddle-focus homoclinics in symmetric systems}
%\title{Intricate order of Shilnikov saddle-focus homoclinics in symmetric systems}

\author{Tingli Xing}
\affiliation{Department of Mathematics \& Statistics, Georgia State University, Atlanta, Georgia 30303, USA}. 
\email{tinglixing107@gmail.com}

\author{Krishna Pusuluri}
\affiliation{Department of Biology, Emory University, Atlanta, Georgia 30322, USA}
\affiliation{Neuroscience Institute, Georgia State University, Atlanta, Georgia 30303, USA}
\email{pusuluri.krishna@gmail.com; Corresponding author}

\author{Andrey L. Shilnikov}
\affiliation{Neuroscience Institute, and Department of Mathematics \& Statistics, Georgia State University, Atlanta, Georgia 30303, USA}
\email{ashilnikov@gsu.edu}

\date{\today}

\begin{abstract}
Using the technique of Poincar\'{e} return maps, we disclose an intricate order of the subsequent homoclinics near the primary homoclinic bifurcation of the Shilnikov saddle-focus in systems with reflection symmetry. We also reveal the admissible shapes of the corresponding bifurcation curves in a parameter plane of such systems. The scalability ratio of geometry and organization is proven to be universal for such homoclinic bifurcations of higher orders. Two applications with similar dynamics due to the Shilnikov saddle-foci, a smooth adaptation of the  Chua circuit and a 3D normal form, are used to illustrate the theory. 
\end{abstract}

%\begin{keywords}
%kneading invariant, symbolic dynamics, T-points, chaos, homoclinic and heteroclinic orbits
%\end{keywords}

\maketitle

{\bf The bifurcation of the Shilnikov saddle-focus is the key for understanding the origin and structure of deterministic chaos in diverse systems including diverse applications from (astro)physics, neuroscience, economics.  
    This article is meant to deepen our understanding of the fine organization of bifurcation unfoldings, including multiple shapes of bifurcation curves in a parameter plane of typical ${\mathbb Z}_2$-symmetric systems. We further develop and showcase the new symbolic approach that lets us disclose a stunning array of homoclinic and heteroclinic bifurcations of the Shilnikov saddle-foci in two representative examples.}      

\section{Introduction}\label{sec:1}

  The aim of this paper is two-fold: its first goal, following the pioneering work of L.P.~Shilnikov on the saddle-focus \cite{Shilnikov1965,LP67,LP68,LP70} and the two later papers~\cite{GTGN97,GLP2007} on its bifurcations, we will begin with examining the structure(s) of homoclinic bifurcation unfoldings in a parameter plane. The second goal is to illustrate computationally the universality and the wealth of such homoclinic bifurcations of the Shilnikov saddle-focus in two representative ODE systems.      
In its second part, this paper is partially an extension of our previous works~\cite{BSS12,Barrio2013,XBS14,pusuluri2017unraveling} on the so-called Lorenz-like systems~\cite{LO63,LP80,ALS86,ASHIL93,SST93} to introduce and demonstrate a new computational  approach~\cite{pusuluri2018homoclinic,pusuluri2019symbolic,Pusuluri2020Chapter,pusulurihomoclinicCNSNS} capitalizing on the symbolic description of homoclinic chaos due to Shilnikov saddle-foci in symmetric systems. An important feature of Lorenz-like systems with partial ${\mathbb Z}_2$-symmetry, i.e., $(x,y,z) \leftrightarrow (-x,-y,z)$, is the universality of complex unfoldings in the parameter space, which are due to the abundance  of homoclinic bifurcations of the plain saddle with a pair of 1D unstable separatrices at the origin. These unfolding are also stirred by the highly characteristic codimension-two T-points, corresponding to the homoclinic connections between the saddle and a pair of symmetric saddle-foci. No 3D Lorenz-like system, except those possessing the full reflection symmetry $(x,y,z) \leftrightarrow (-x,-y,-z)$, can accommodate a saddle-focus with 1D outgoing separatrices and a 2D stable manifold that is due specifically to a pair of complex conjugate characteristic exponents, with a negative real part.  It is well-known that the occurrence of a single homoclinic orbit of the Shilnikov saddle-focus can give rise to the onset of chaotic dynamics, including countably many nearby periodic orbits in the phase space of such a system. Shilnikov's theory from the 60's demonstrated the significance of the organizing role of homoclinic orbits in the hierarchy of deterministic chaos 
  \cite{AfrShil1983}.

\begin{figure*}[thb!]
\begin{center}
\includegraphics[width=.9\linewidth]{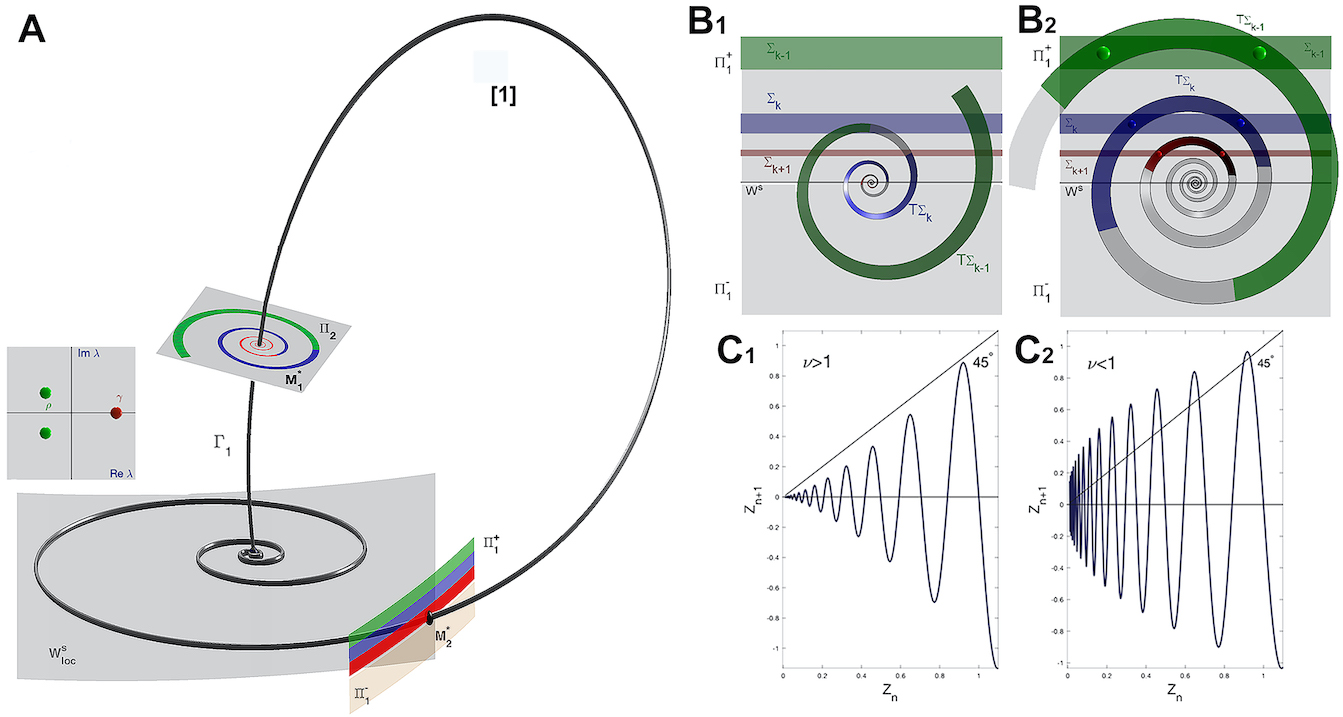}
\caption{(A) Primary homoclinic orbit $ \bar \Gamma_1$ of a saddle-focus of (2,1)-type, i.e., with 2D stable manifold $W^s$ and 1D unstable manifold $W^u$ in ${\mathbb R}^3$. The colored stripes on a local 2D cross-section $\Pi_1$ transverse to  $W^s$ are taken along the trajectories near the equilibrium state to be transformed into a spiral on the local cross-section $\Pi_2$ transverse to $W^u$. 
 (B) 2D Poincar\'e return map $\Pi_2 \to \Pi_1$ is a contraction in (B$_1$) with a saddle index $\nu >1$ (corresponding to 1D map shown in (C$_1$) ), or   an expansion in (B$_2$) with non-empty intersections $T\, \Sigma_k \cap \Sigma_k $ giving rise to the onset of countably many Smale horseshoes and saddle POs corresponding to repelling fixed points (FP) in the 1D map in (C$_2$) when $\nu <1$    -- the so-called Shilnikov condition. From \url{http://www.scholarpedia.org/article/Shilnikov_bifurcation}.}\label{fig1}
\end{center}
\end{figure*}

Let us re-iterate without excessive detail, what is well-known about the Shilnikov homoclinic saddle
focus. The Reader is welcome to consult with L.P.~Shilnikov's original papers~\cite{Shilnikov1965, LP67,LP68,LP69,Shilnikov_scolarpedia,Shilnikov_heritage} and his co-authored textbooks~\cite{a2001methods,arnold2013dynamical}, as well as with other relevant papers on the theory \cite{G83,Bel84,OSh86e,ovsyannikov1992systems,GTGN97,GLP2007}  by his students,  and its various extensions~\cite{ArnCoulTres1981,gn83,gkn84,medrano2005basic,fowler1991bifocal} and diverse applications  \cite{ArgArnEleRich1993,B00,FeudelPei2000, nicolay2004low, BarrioShilnikov2011,KopGasSlu1992,BasHud1988, toniolo2005phenomenological, cortes2013short,Rossler2020,barnett2020shilnikov}. Figure~\ref{fig1}A demonstrates the simplest or primary homoclinic orbit to a saddle-focus of the topological (2,1)-type; more details including analytical results will be given in Section 2. Here, the topological (2,1)-type means that the saddle-focus has a pair of complex conjugate characteristic exponents
(small green dots in the inset of Fig.~\ref{fig1}) in the left open complex half-plane, and one positive real one. 
To be the Shilnikov saddle-focus, the complex pair is to be the closest to the imaginary axis.  Trajectories of a system near such a  saddle-focus take a local cross-section $\Pi_1^+$ transverse to the 2D stable manifold $W^s_{loc}$ and map onto another cross-section $\Pi_2$ transverse to a 1D unstable separatrix $\Gamma_1$. Then, the colored stripes on $\Pi_1^+$ will be transformed into a spiral, sometimes called the Shilnikov snake with an ordered color pattern, on $\Pi_2$. Next, the global map takes the spiral and maps it back onto the first cross-section as shown in Fig.~\ref{fig1}B. Depending on the ratio of the local stability to instability at the saddle-focus, there are two options. One is when stability exceeds instability, the overall map is a contraction (Fig.~\ref{fig1}B$_1$); otherwise, it is an expansion, see Fig.~\ref{fig1}B$_2$. The latter implies that the colored (green, blue and red) stripes in $\Pi_1^+$ can be reached and crossed, geometrically, by the their arched images on the spiral. Such crossings are interpreted as the formation of countably many topological Smale horseshoes,  giving rise to countably many unstable periodic orbits, and the onset of complex shift dynamics just near the primary homoclinic orbit in the phase space of the given system. The corresponding 1D return maps are shown in Fig.~\ref{fig1}C. These are basically the ``parameterizations'' of the spirals on either coordinate axes. One can see from Fig.~\ref{fig1}C$_1$ that the contraction map, when shifted up, will produce a single stable fixed point (FP) at the intersection with the $45^\circ$-line from the origin in the 1D return map, which corresponds to the saddle-focus in the phase space. On the contrary, the expansion map in Fig.~\ref{fig1}C$_2$ with characteristic oscillations generates countably many crossings, read FPs, on the $45^\circ$-line. When the homoclinic orbit in Fig.~\ref{fig1}A splits above/below $W^s_{loc}$, the 1D return map shifted up/down perseveres most of the FPs. Some of its oscillations will become tangent to the $45^\circ$-line to produce new crossings. Such tangencies cause saddle-node bifurcations, soon to be followed by period-doubling ones.  This is a reason why the Shilnikov bifurcation in systems with 3D phase space is a precursor of deterministic chaos, associated with the so-called quasi-chaotic attractors in which hyperbolic subsets coexist with stable periodic orbits emerging through saddle-node bifurcations.

\begin{figure}[t!]
\begin{center}
\includegraphics[width=.99\linewidth]{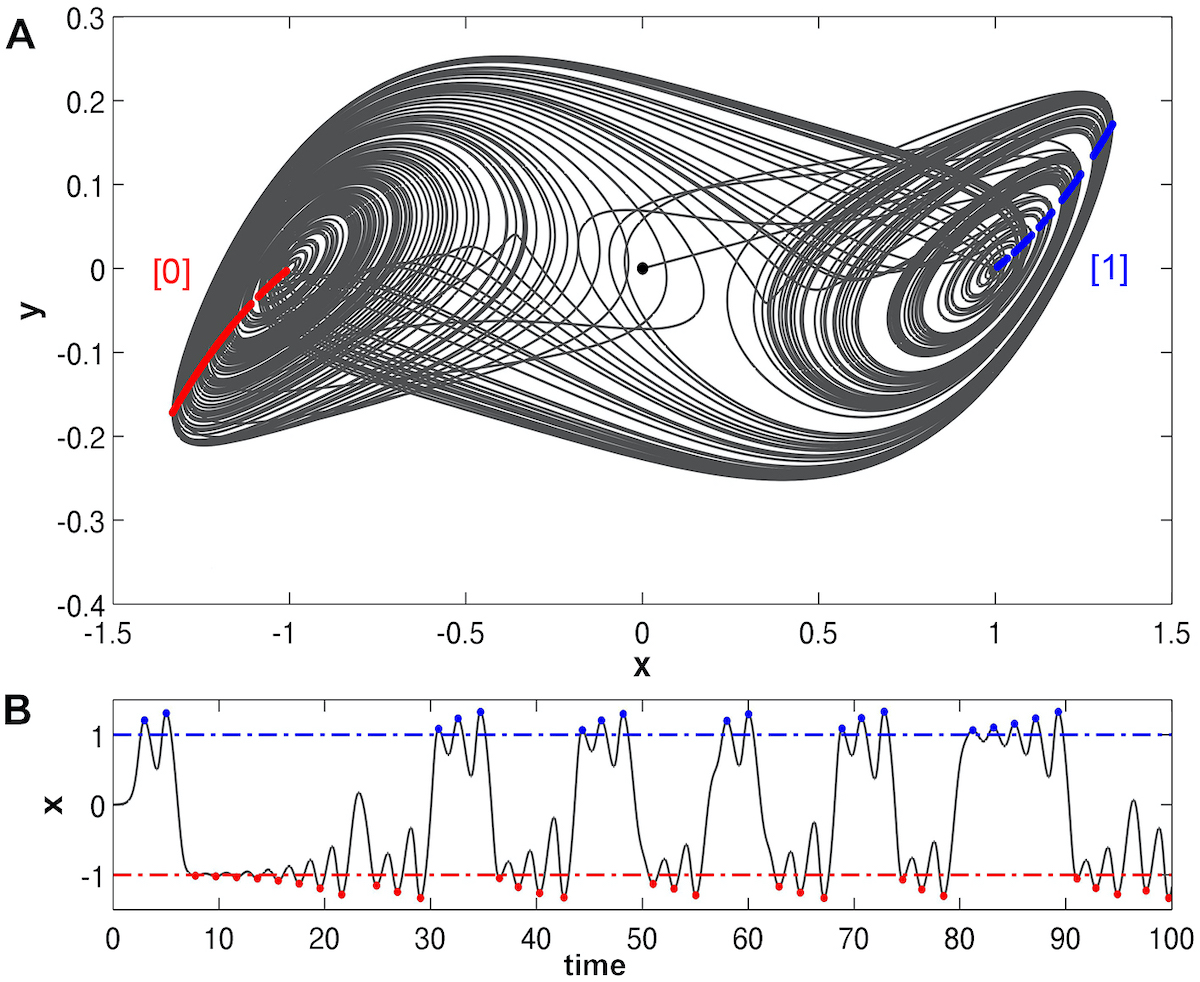}
\caption{(A) Homoclinic chaos due to three Shilnikov saddle-foci in the symmetric 3D Chua system~(1) at $a=10.16$ and $b=14.7$.  The right 1D separatrix of the saddle-focus $O$ at the origin (black dot) fills in the double-scroll  attractor by making an unpredictable number of turns around and by switching between two other saddle-foci $O_{1,2}(\pm 1, 0, 0)$, separated by our key player -- the saddle-focus of the type (2,1) at the origin in the 3D phase space. Its symbolic, binary sequence is generated using a simple rule: [1] or [0] whenever the trajectory turns around $O_1$ or $O_2$, respectively. Or alternatively, when its $x$-coordinate reaches a next successive maximum/minimum above/below +1/-1, respectively. See the progression of $x(t)$ in (B).}\label{fig2}
 \end{center}
\end{figure}

A representative example of deterministic chaos in due to three Shilnikov saddle-foci in the phase space (of the smooth Chua model below) is depicted in Fig.~\ref{fig2}. This figure also illustrates the concept of  $\{0,\,1\}$-based binary symbolic description in application to symmetric systems with chaotic dynamics.

 As pointed our earlier, this paper includes two parts: a theoretical one followed by computational sections.  First, we extend the theory to analytically disclose the structure of local bifurcation unfolding of subsequent homoclinic bifurcations of the Shilnikov saddle-focus, near the primary one, see Figs.~\ref{fig1} and \ref{fig5}, in ${\mathbb Z}_2$-symmetric systems. The second goal is to reveal the global fine organization of chaos due to  the Shilnikov saddle-focus homoclinic bifurcations in two exemplary, ${\mathbb Z}_2$-symmetric systems, through detailed visualizations   with the aid of a newly proposed computational approach capitalizing on the symbolic description of trajectories on observable strange attractors.

\begin{figure*}[t!]
\begin{center}
\includegraphics[width=.8\linewidth]{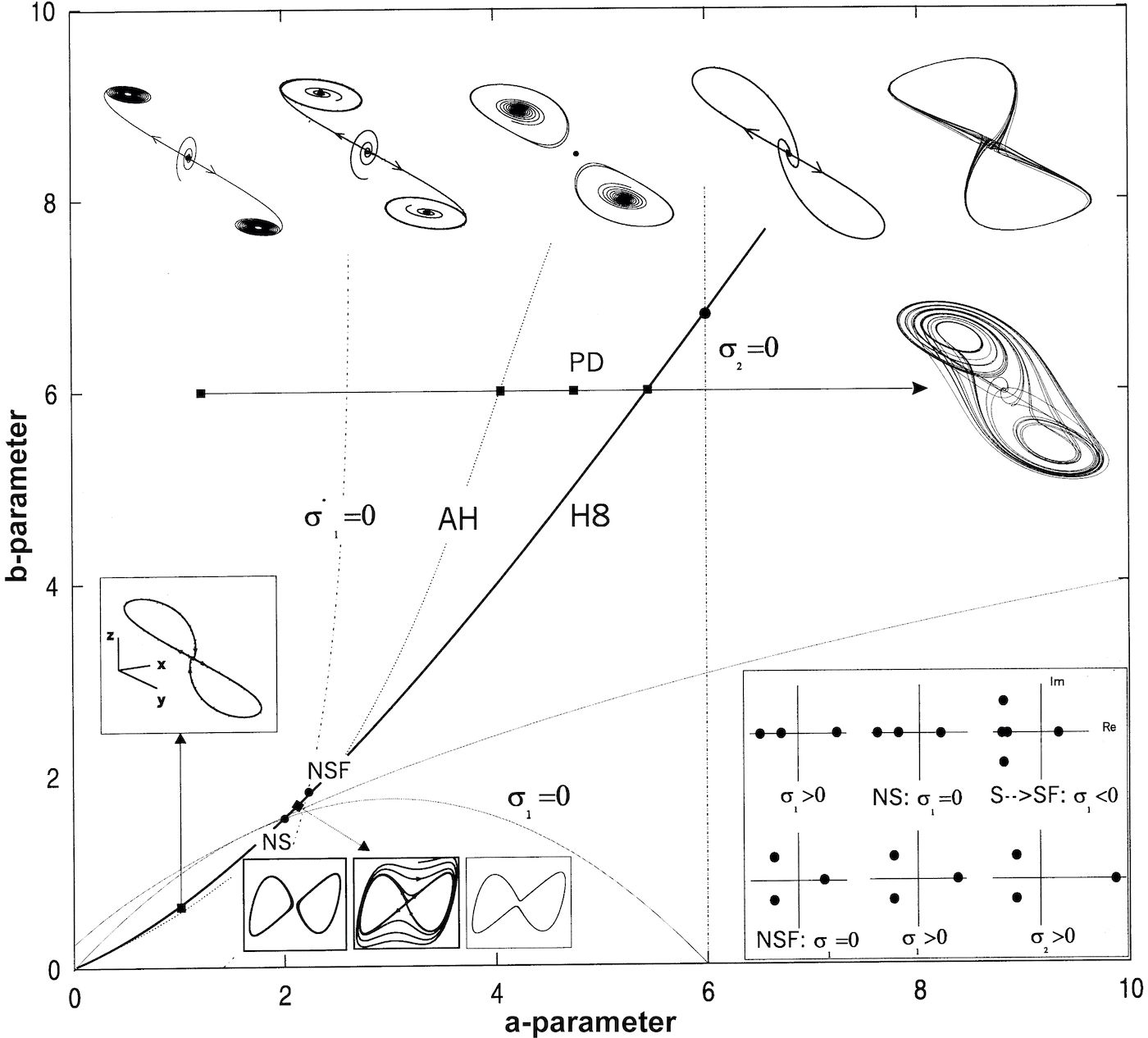}
\caption{Bifurcation diagram of the cubic Chua model. Abbreviations AH and H8 stand for, resp., a supercritical Andronov-Hopf bifurcation of the symmetric equilibria $O_{1,2}$ and a figure-8 homoclinic connection of the saddle-focus $O$ at the origin, while NS, NSF, and $S \rightarrow SF$ stand for a  neutral (resonant) saddle-focus with a zero saddle value ($\sigma_1=0$), a saddle-focus with zero divergence ($\sigma_2=0$) and saddle $\Leftrightarrow $ saddle-focus transition (see the panel at the right-bottom corner), resp. Along the pathway $b=6$, firstly, the 1D unstable separatrices $\Gamma_{1,2}$ of the origin $O$ converge to stable equilibria $O_{1,2}$; next, as the $a$-parameter is increased, they converge to stable periodic orbits (POs) that further lose stability through a first period-doubling (PD) bifurcation in the forthcoming cascade.  These saddle POs become the primary homoclinic loops of the Shilnikov saddle-focus $O$, thus giving rise to the onset of homoclinic chaos in the model~(\ref{ch1}); from Ref.~\cite{a2001methods}}\label{fig3}
\end{center}
\end{figure*}

The first example is a smooth approximation of the Chua's circuit~\cite{MTC84}.  The circuit, including two capacitors, two resistors, one inductor, and a nonlinear element, Chua's diode, is described by a 3D system of ODEs, with a a single nonlinear term.  All of its parameters have specific physical meanings~\cite{RG10}. Originally, the non-linearity was described using a piece-wise function, that was later replaced with a smooth cubic function in Ref.~\cite{T05} Both  systems were compared in detail in Ref.~\cite{RG10}

The smooth Chua model with a cubic nonlinearity~\cite{KRC93,B98} is given by  
\begin{equation}
\dot{x} = a\left (y+\frac{x}{6}-\frac{x^3}{6} \right ), \quad \dot{y}=x-y+z, \quad \dot{z}=-by,   \label{ch1}
\end{equation}
with $a,\,b>0$ being bifurcation parameters. The system is reflection or  ${\mathbb Z}_2$--symmetric, i.e., invariant under the involution $(x, y, z) \rightarrow (-x, -y, -z)$.  It has three equilibrium states: $O(0,0,0)$ can  be a saddle of the topological type (2,1), i.e. with two 1D unstable separatrices, call them $\Gamma_1$ and $\Gamma_2$ and a 2D stable manifold $W^s$, or a saddle-focus of the same topological type, while $O_1(-1,0,1)$ and $O_2(1,0,-1)$ can be stable or saddle-foci of the type (1,2). In the chaotic region of our particular interest in the parameter plane, all three equilibria are saddle-foci. Figure~\ref{fig3} illustrates a bifurcation diagram for the equilibrium states in the cubic Chua model~(\ref{ch1}), see Refs.~\cite{KRC93,a2001methods} for more details. 
 
The other example employed for the illustration of our symbolic approach to disclose the global organization of homoclinic and heteroclinic bifurcations of the Shilnikov saddle-foci is an asymptotic normal form~\cite{CTA1979}:  
\begin{equation}\label{act}
\dot{x} = y, \quad \dot{y}=z, \quad \dot{z}=-b\,z -y + a\,x \left (1-x^2 \right ) 
\end{equation}
with $(a,\,b)>0$ being the bifurcation parameters, describing a local bifurcation unfolding in systems, near an equilibrium state with a triplet of zero characteristic exponents on a  ${\mathbb Z}_2$-symmetric central manifold. Its phase space with three saddle-foci may look similar to that of the cubic Chua model~(\ref{ch1}).  This normal form, as well as some other  systems were in-detail studied in Ref.~\cite{ArnCoulTres1980,ArnCoulTres1981,ArnCoulTres1982,ArnCoulSpTres1985}, which  along with the Brussel group~\cite{G83,gn83,gkn84}, were the very first works in the West that began studying the Shilnikov saddle-focus and spiral chaos around it.  We will refer to Eqs.~(\ref{act}) as the cubic Arneodo-Coullet-Spiegel-Tresser (ACST) model after the authors of the series of the publications.

\begin{figure}[hbt!]
\begin{center}
\includegraphics[width=.999\linewidth]{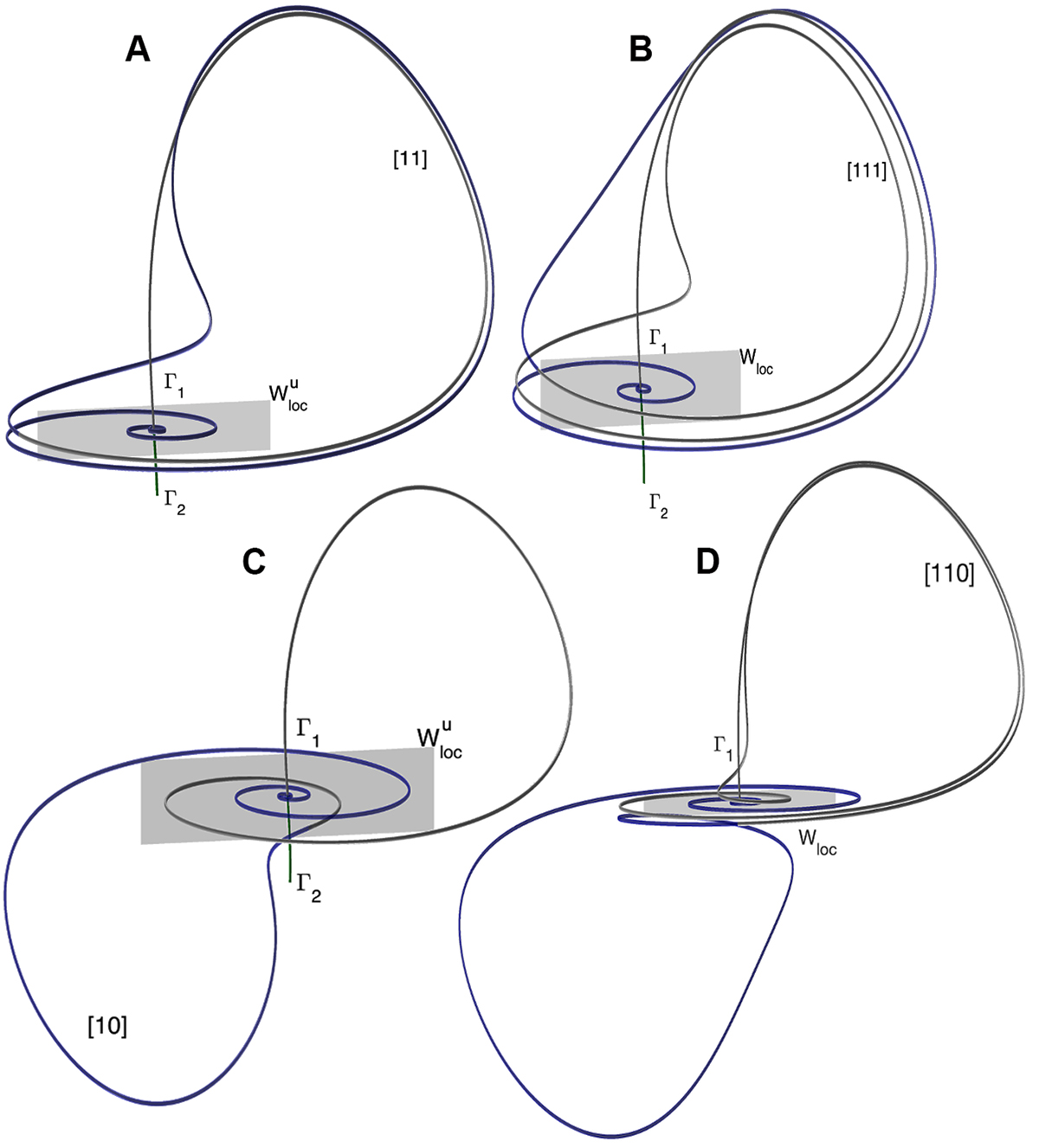}
\caption{Secondary homoclinic orbits of a saddle-focus~$(2,\,1)$ in the phase space of a ${\mathbb Z}_2$-symmetric system:  one-sided double (A) and triple (B) homoclinics encoded symbolically  as $[11]$ and $[111]$, resp. (C,\,D) Figure-8 homoclinic orbits (of $\Gamma_1$) with encoded as $[10]$ and $[110]$, resp.  The depicted homoclinic orbits were generated, for sake of the illustration purpose, by solutions of a 4D ${\mathbb Z}_2$-symmetric Shimizu-Morioka type model with the Shilnikov saddle-focus at the origin, see Eqs.~(C.7.7) on page 924 in Ref.~\cite{a2001methods}}\label{fig4}
\end{center}
\end{figure}

Unlike hyperbolic systems with complex dynamics, the Chua model constantly  undergoes abrupt transitions or
bifurcations, including homoclinic, saddle-node and period-doubling, within a parameter region of the existence of the three Shilnikov saddle-foci, see Ref.\cite{GTGN97} However, we will not discuss the one sided chaos due to the two symmetric saddle-foci, because its bifurcation structure is essentially the same
as observed in the R\"ossler system, see Refs.~\cite{Rossler1976,Rossler2020} Instead, we focus on the symmetric and asymmetric homoclinics generated by the central saddle-focus at the origin, and
how their structures are embedded in the parameter space. We use binary symbols $[0]$ and $[1]$ to symbolically encode
such homoclinic orbits, see Fig.~\ref{fig4} illustrating the concept.  
 Namely,  [1] is used to describe the passes of the separatrix $\Gamma_1$, while [0] is reserved for the other symmetric separatrix  $\Gamma_2$. For example, double or triple {\em  one-sided} homoclinic orbits  are encoded as [11] or [111], or symmetrically as [00] or [000], resp., see  Fig.~\ref{fig4}A-B. If $\Gamma_1$ misses the primary loop and goes underneath the stable manifold $W^s_{loc}$  before it comes back to the saddle-focus as illustrated in Fig.~\ref{fig4}C, then its code is [10]. Figure~\ref{fig4}D depicts a more complex triple homoclinic orbit encoded as [110].    

The paper is organized as follows. In Section 2, we will present our analytical results on homoclinic bifurcations of the Shilnikov saddle-focus in reflection-symmetric systems. Section 3 will introduce a symbolic computational tool (see also Refs.\cite{BSS12,XBS14}) and apply it to the smooth Chua model~(\ref{ch1}) to compare numerical findings with the theoretical results from Section~2. Section 3 will focus on the numerical study of bi-parametric sweeps of the normal form~(\ref{act}), which is followed by conclusions and discussion.

\section{Analytical approach: homoclinic bifurcation structure of the Shilnikov saddle-focus in symmetric systems}\label{sec:2}

\begin{figure}[hbt!]
\begin{center}
\includegraphics[width=.85\linewidth]{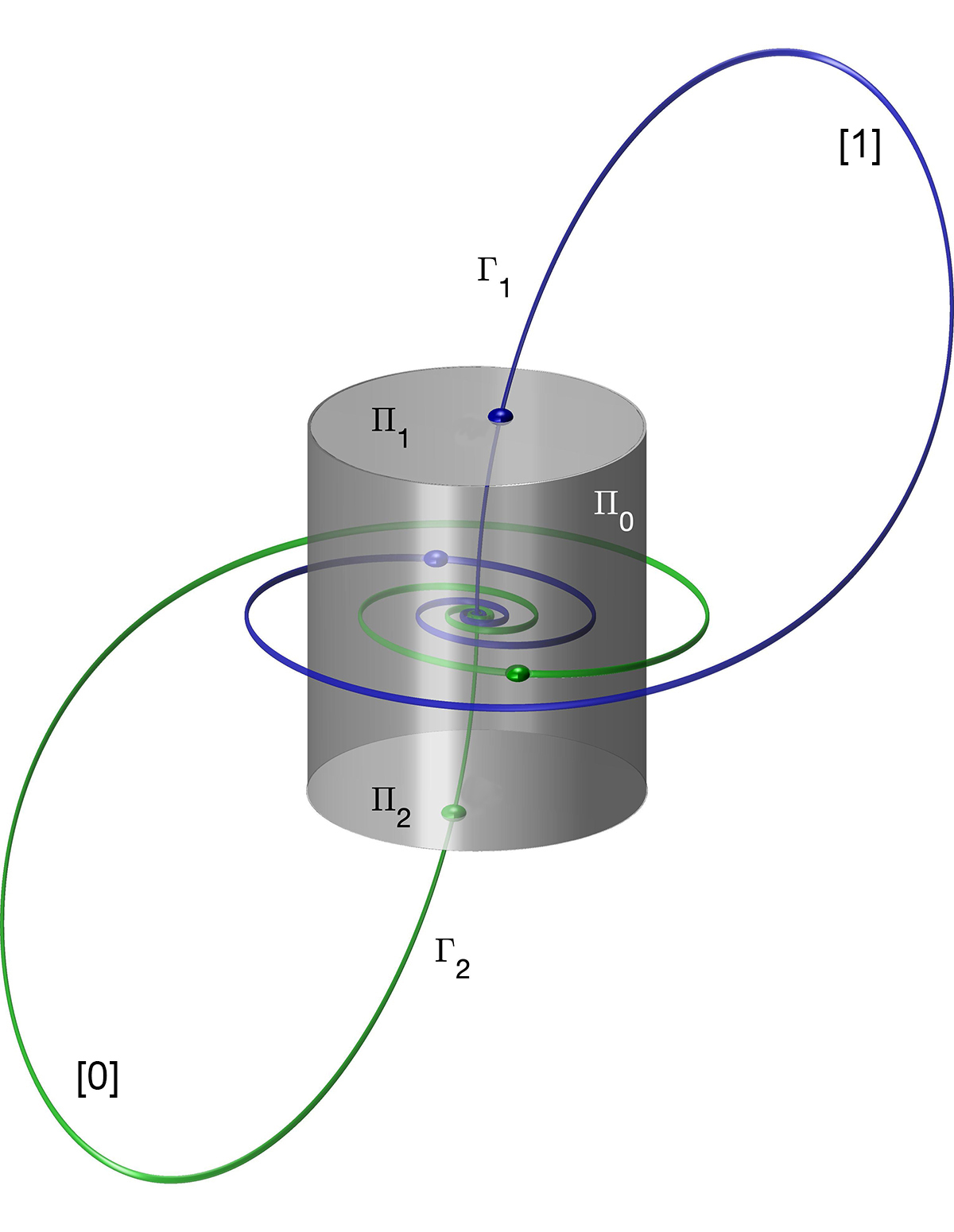}
\caption{A pair of the primary homoclinic orbits, symbolically encoded as [1] and [0], of a saddle-focus at the origin in the phase space of a ${\mathbb Z}_2$-symmetric system. The 1D outgoing separatrices $\Gamma_1$ and $\Gamma_2$ hit outwardly the top $\Pi_1$ and bottom $\Pi_2$ bases of a cylinder-shaped cross-section enclosing the origin and its sidewall $\Pi_0$ upon their return to the saddle-focus.}\label{fig5}
\end{center}
\end{figure}

Let us consider the homoclinic Shilnikov saddle-focus of the (2,1)-type at the origin $O$ of a 3D system with reflection symmetry. Figure~\ref{fig5} illustrates this where both 1D separatrices $\Gamma_{1,2}$ leave the saddle-focus symmetrically, and after a short excursion, come back to it along its 2D stable manifold $W_O^s$. This is called a primary homoclinic figure-8.  In what follows, we will consider how small smooth perturbations of a system with such a figure-8 can generate longer subsequent homoclinic orbits of the saddle-focus, under the fulfillment of a single so-called Shilnikov condition~\cite{LP68}. We will also describe how such homoclinic bifurcations are embedded in a parametric plane.

Following Ref.\cite{LP68,a2001methods}, let us use the following form of a $\mathbb{Z}_2$-symmetric system near the saddle-focus: 
\begin{equation}
\begin{array}{lcl}
\dot{x} &=& -\rho(\mu) x-\omega(\mu)y+F_1(x,y,z,\mu), \\
\dot{y} &=& ~~\omega(\mu)x-\rho(\mu)y+F_2(x,y,z,\mu), \\
\dot{z} &=& ~~\lambda(\mu)z+F_3(x,y,z,\mu),
\end{array} 
\label{ch3}
\end{equation}
where $F_i$ are smooth functions, so that $F_i(\cdot)=-F_i(-(\cdot))$ and they and their first derivatives vanish at $O(0,0,0)$ for all small $\mu$;  the primary homoclinic figure-8 occurs at $\mu=0$. The characteristic exponents of the saddle-focus are given by $-\rho(\mu) \pm i \omega(\mu)$ so that  $\rho(\mu)$ and  $\omega(\mu)>0$, and $\lambda(\mu)>0$. The so-called {\em saddle index} is given by $\nu(\mu)=\frac{\rho(\mu)}{\lambda(\mu)}<1$; this is the Shilnikov condition~\cite{LP68} needed for complex dynamics of the finite-shift type to merge in a system with such a saddle-focus. In this normalized system, the $z$-axis is the linearized unstable manifold $W_O^u$ and the $(x,y)$-plane is the linearized stable manifold $W_O^s$ of the saddle-focus at the origin. The solution of the linearized system (\ref{ch3}) initiated at a point $(x_0, y_0, z_0)$ can be written as
\begin{equation}
\begin{array}{ccl}
x(t) &=& e^{-\rho(\mu)t}\,[x_0\cos(\omega(\mu) t)-y_0\sin(\omega(\mu) t)], \\
y(t) &=& e^{-\rho(\mu)t}\,[y_0\cos(\omega(\mu) t)+x_0\sin(\omega(\mu) t)], \\
z(t) &=& e^{\lambda(\mu)t}\,z_0.
\end{array}
\label{ch4}
\end{equation}
As a transverse cross-section, we pick a sufficiently small cylinder (see Fig.~\ref{fig5}) enclosing the saddle-focus, to construct a Poincar\'{e} return map in cylinder coordinates  $(r, \varphi, z)$, following Ref.~\cite{G83} It is constructed with a side wall given by  $\Pi_0: r=R (-R<z<R)$, and top and bottom disks given by $\Pi_1: z=R (0<r<R)$ and $\Pi_2: z=-R (0<r<R)$, respectively; here, $R$ is sufficiently small. 

For $z_0>0$, the local map $T_0: \Pi_0 \mapsto \Pi_1((\varphi_0, z_0) \mapsto (r,\theta))$ is calculated from Eq.~(\ref{ch4}), noting that $x_0=R\cos\varphi_0$, $y_0=R\sin\varphi_0$,  $x(t)=r\cos\theta$, $y(t)=r\sin\theta$ and $z(t)=R$. It is given by 
\begin{equation} 
T_0: \quad \left [ 
\begin{array}{lcl}
r &=& R(z_0/R)^{\nu(\mu)}, \\ 
\theta &=& \varphi_0+(\omega(\mu)/\lambda(\mu) )\ln(R/z_0). 
\end{array} \right .
\label{ch5}
\end{equation}

 Similarly, when $z_0<0$, the local map $T_0': \Pi_0 \mapsto \Pi_2((\varphi_0, z_0) \mapsto (r,\theta))$ can be calculated from (\ref{ch4}) as:
$$T_0': \left [ 
\begin{array}{lcl}
 r &=& R(-z_0/R)^{\nu(\mu)}, \\
 \theta &=& \varphi_0+(\omega(\mu)/\lambda(\mu) )\ln(-R/z_0). 
\end{array} \right .
$$
The global map $T_1:~\Pi_1 \mapsto \Pi_0$  (which is $(x,y) \mapsto (\varphi_0, z_0)$ or $(r\cos\theta, r\sin\theta) \mapsto (\varphi_0,z_0)$) along the separatrices $\Gamma_{1,2}$, returning to the cylinder-shaped cross-section, can be approximated by a linear transformation:
\begin{equation} T_1:~\left [ 
\begin{array}{lcl}
\varphi_0 &= &a_1\mu+a(\mu)x+b(\mu)y,\\
 &=& a_1\mu+A(\mu)r\cos(\theta+\alpha_1(\mu))+O(r^2),\\
z_0 &=& \mu+c(\mu)x+d(\mu)y,\\
 &=& \mu+B(\mu)r\sin(\theta+\alpha_2(\mu))+O(r^2),
\end{array} 
\right .
\label{ch6}\end{equation}
where $A(0)B(0)\cos[\alpha_1(0)-\alpha_2(0)] \neq 0$ for a non-degenerate linear transformation. The map $T_1': \Pi_2 \mapsto \Pi_0$ can be derived from $T_1$ using reflection symmetry. For $(x,y) \in \Pi_2$, $T_1'(x,y)$ and $T_1(-x,-y)$ are symmetric with respect to the origin, therefore:
$$T_1': \quad \left [ 
\begin{array}{rcl}
\varphi_0 &=& \pi+a_1\mu-A(\mu)r\cos(\theta+\alpha_1(\mu))+O(r^2),\\
    z_0  &=& -\mu+B(\mu)r\sin(\theta+\alpha_2(\mu))+O(r^2).
\end{array} \right .
$$
\begin{figure*}[hbt!]
\begin{center}
\includegraphics[width=.85\linewidth]{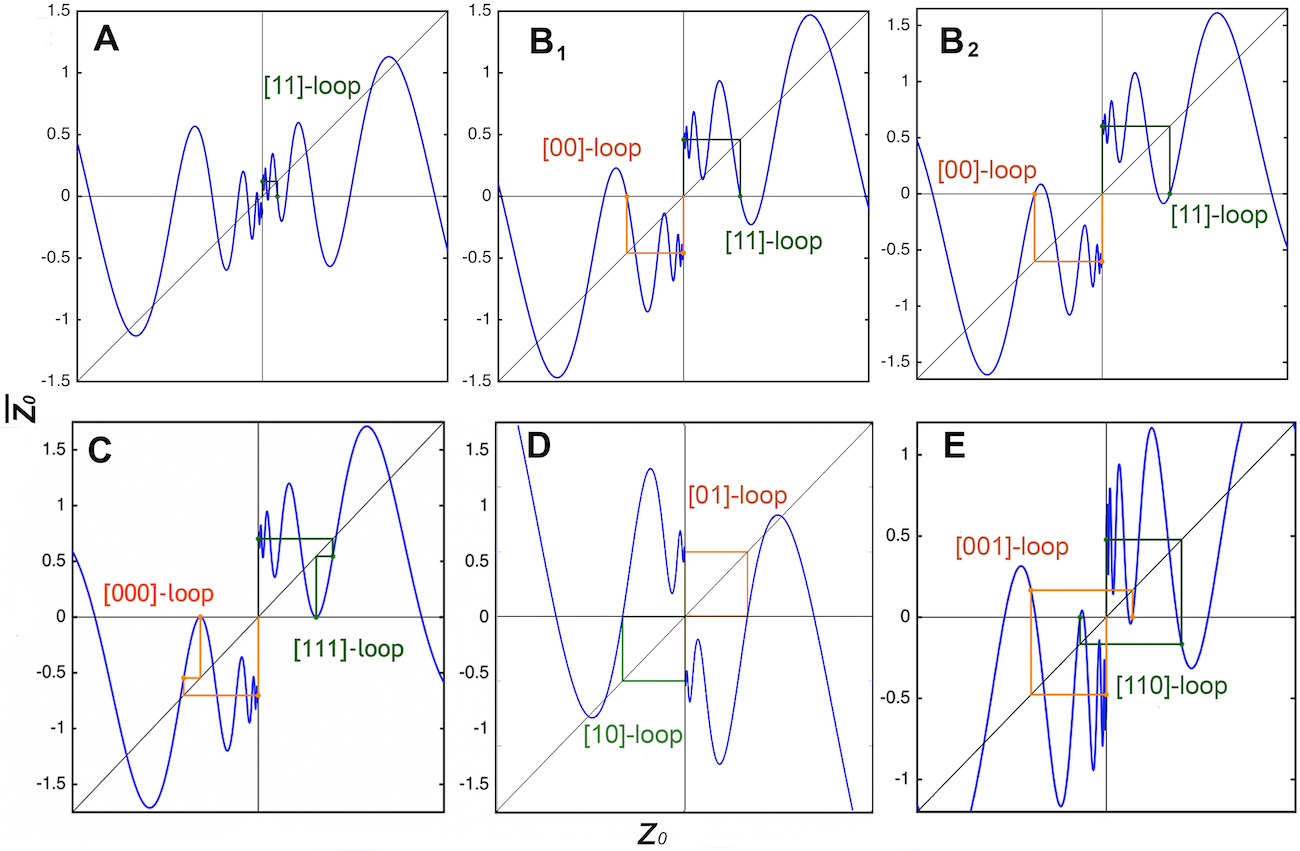}
\caption{1D Poincar\'{e} return maps associated with the subsequent homoclinic orbits of the saddle-focus in the phase space shown in Fig.~\ref{fig4}: the forward iterates of the origin correspond to the 1D unstable separatrices returning to the saddle-focus. (A)-(B) Multiplicity of symmetric double homoclinic orbits encoded with [11] and  [00] are due to countably many zeros in the 1D oscillatory map near such a saddle-focus. (C) Triple homoclinic orbits encoded with [111] and [000], like the homoclinics shown in Fig.~\ref{fig4}B. (D)-(E) Unlike {\em one-sided} orbits,  these 1D maps generate a multiplicity of more complex pairs of figure-8 orbits such as [10]/[01] and [110]/[001], similar to those shown in Figs.~\ref{fig4}C and D, resp.}\label{fig6}
\end{center}
\end{figure*}
Hence, when $z_0>0$, the complete return map $T=T_1 \circ T_0 : \Pi_0 \mapsto \Pi_0$ is given by:
$$ T: \quad \left [  
\begin{array}{lcl}
\bar{\varphi_0}&=& a_1\mu+A(\mu)R(z_0/R)^{\nu(\mu)}\cos(\varphi_0+\frac{\omega(\mu)}{\lambda(\mu) }\ln(R/z_0) \\ &&+\alpha_1(\mu))+O(z_0^{2\nu_(\mu)}),\\
\overline{z_0}&=&\mu+B(\mu)R(z_0/R)^{\nu(\mu)}\sin(\varphi_0+\frac{\omega(\mu)}{\lambda(\mu) }\ln(R/z_0)\\&&+\alpha_2(\mu))+O(z_0^{2\nu(\mu)}).
\end{array}  \right .
$$

For the case $z_0<0$, the corresponding return map $T'=T_1' \circ T_0' : \Pi_0 \mapsto \Pi_0$ is given by: 
$$ T': \quad \left [
\begin{array}{lcl}
\tilde{\varphi_0}&=&\pi+a_1\mu-A(\mu)R(-z_0/R)^{\nu(\mu)}\cos(\varphi_0 \\&&+\frac{\omega(\mu)}{\lambda(\mu) }\ln(-R/z_0)+\alpha_1(\mu))+O(z_0^{2\nu_(\mu)}),\\
\bar{z_0}&=&-\mu+B(\mu)R(-z_0/R)^{\nu(\mu)}\sin(\varphi_0 \\&&+
 \frac{\omega(\mu)}{\lambda(\mu) }\ln(-R/z_0)+\alpha_2(\mu))+O(z_0^{2\nu(\mu)}).
\end{array} \right .
$$

Let $A_0=A(0)$, $B_0=B(0)$, $\Omega_0=\omega(0)/\lambda(0) $,  $\nu_0=\nu(0)$, $\phi_1=-\alpha_1(0)-\Omega_0\ln R$ and $\phi_2=-\alpha_2(0)-\Omega_0\ln R$. Keeping only the dominant terms, these maps can be simplified as follows:
\begin{equation} T: \left [ 
\begin{array}{lcl}
\bar{\varphi_0} &= &a_1\mu+A_0R^{1-\nu_0} z_0^{\nu_0}\cos(\Omega_0\ln{z_0}+\phi_1-\varphi_0)+ \\ &&+O(z_0^{2\nu_0}),\\
\bar{z_0}&=&\mu-B_0R^{1-\nu_0} z_0^{\nu_0}\sin(\Omega_0\ln{z_0}+\phi_2-\varphi_0)+\\&&+ O(z_0^{2\nu_0}),
\end{array} \right .
\label{ch7}
\end{equation}
and
\begin{equation} T': \left [ 
\begin{array}{lcl}
\tilde{\varphi_0}&=&\pi+a_1\mu-A_0R^{1-\nu_0}(-z_0)^{\nu_0}\cos(\Omega_0\ln (-z_0)+\\ && +~ \phi_1-\varphi_0)+O(z_0^{2\nu_0}),\\
\tilde{z_0}&=&-\mu-B_0R^{1-\nu_0} (-z_0)^{\nu_0}\sin(\Omega_0\ln (-z_0)+\\&&+ \phi_2-\varphi_0)+O(z_0^{2\nu_0}).
\end{array} \right . 
\label{ch8}
\end{equation}

A homoclinic orbit that passes $l$ times through the cylinder-wall $\Pi_0$ is called an $l$-loop homoclinic orbit, while shorter ones with $2$ or $3$ passes are called double- or triple-loop homoclinic orbits, respectively.
Whenever either 1D separatrix $\Gamma_i$  hits $\Pi_0$  with $z>0$, we extend its encoding with the symbol $1$;  otherwise, if $z<0$, its code is extended with the symbol $0$.  For example, a homoclinic orbit that passes  through $\Pi_0$ $2$ times with $z>0$, is called a one-sided double separatrix $[11]$-loop/homoclinic orbit (see Fig.~\ref{fig4}A; and ~\ref{fig4}B-D for longer homoclinic orbits and their symbolic codes). The $\mu$-parameter is often referred to as a splitting parameter whose positive/negative variations split the primary homoclinic orbit, say $\bar \Gamma_1$, upward/downward with respect to the saddle-focus or its stable manifold $W^s_{loc}$. 

\begin{figure*}[htb!]
\begin{center}
\includegraphics[width=.85\linewidth]{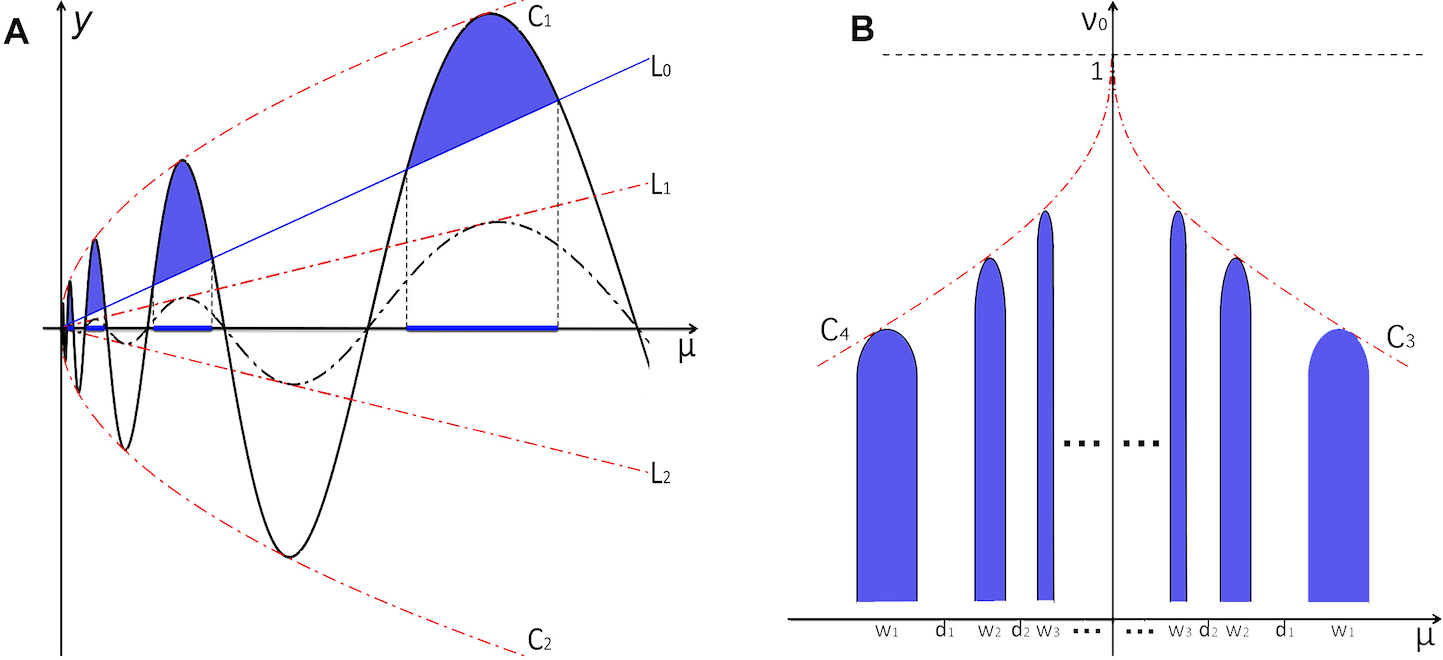}
\caption{Panels illustrating the causality and structure of the bifurcation unfolding of the 1D Poincar\'e return map~(\ref{1dmap}) in the region of $z_1<0$, when $B_0<1$. (A) Condition $z_1<0$ is fulfilled on the graph segments where the sine-function is above the line $L_0$. The union of the projected (blue) intervals yields the range of $\mu$-values resulting in $z_1<0$. The interval endpoints correspond to the formation of the double [11]-and [00]-homoclinic orbits. As $\nu_0 \rightarrow 1$, the amplitude $C_{1,2}$ of the sine-function flattens to two lines, $L_{1,2}$, causing the blue intervals to vanish. (B) Sketch of the $(\mu,\nu_0)$-plane with the colored regions, aka blue $\cap$-bars, in which $z_1<0$. The outlines of the blue $\cap$-bars on the $\mu>0$-side corresponds to the double [11]/[00] homoclinic orbits. The outlines of the blue $\cap$-bars on the $\mu<0$ side correspond to the double [10]/[01]-orbits. The blue bars are bounded from above by the $C_{3,4}$-curves; $w_j$ and $d_j$ stand for the widths and the distances between the bars, resp. }\label{fig7}
\end{center}
\end{figure*}

\begin{figure*}[thb!]
\begin{center}
\includegraphics[width=.85\linewidth]{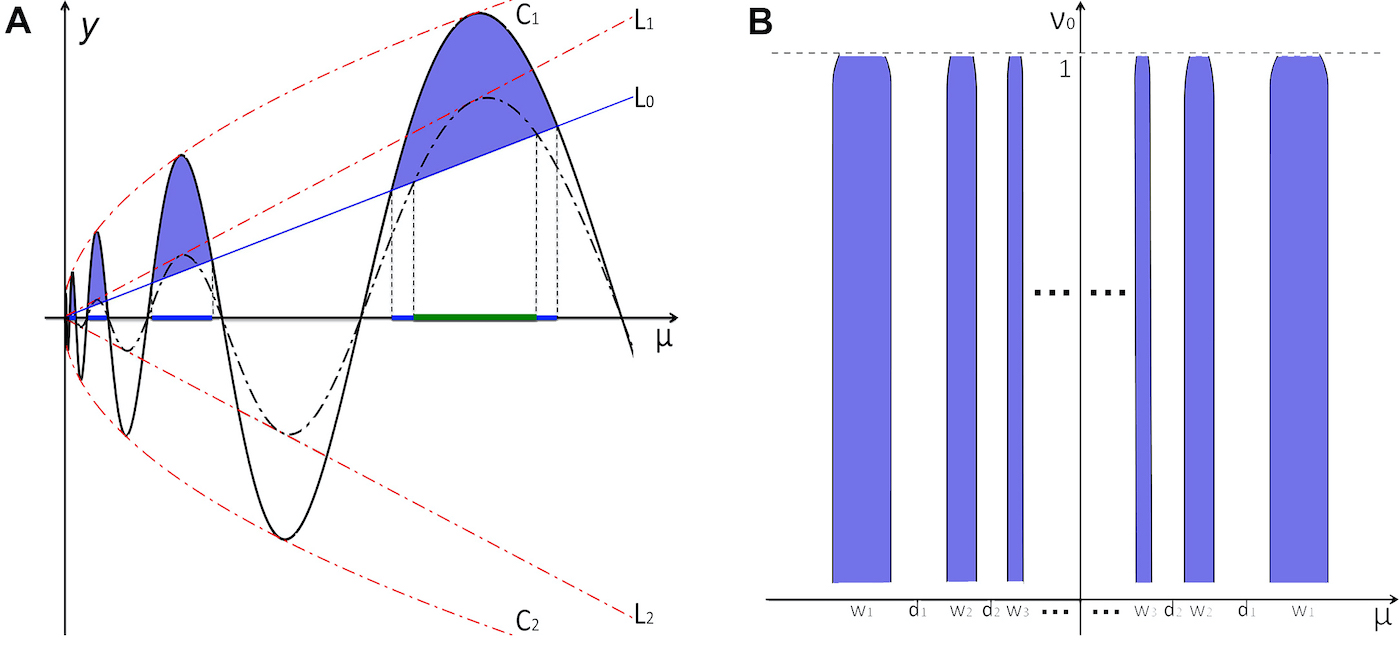}
\caption{Panels illustrating  the causality and structure of the bifurcation unfolding of the 1D Poincar\'e return map~(\ref{1dmap}) in the region of $z_1<0$, when $B_0>1$. (A) $z_1<0$ is fulfilled on the graph segments where the sine-function is above the line $L_0$. Blue intervals on the $\mu$-axis are where $z_1<0$. The interval endpoint are the bifurcation $\mu$-values corresponding to one-sided [11]/[00] homoclinic orbits. As $\nu \rightarrow 1$, the envelopes $C_{1,2}$ of the sine-function flatten to the $L_{1,2}$-lines, thereby making the blue intervals  shorten to the green intervals. (B) The plot sketches the regions  in the $(\mu,\nu_0)$-parameter plane -- blue $\Pi$-bars -- where $z_1<0$.  The outlines of the  blue $\Pi$-bars on the $\mu>0$-side are the bifurcation curves of the double [11]/[00]-orbits. The outlines of the blue  $\Pi$-bars on the $\mu<0$-side correspond to [10]/[10] double homoclinic orbits; here $w_j$ and $d_j$ stand for the widths and the distances between the bars, resp.}\label{fig8}
\end{center}
\end{figure*}

Figure~\ref{fig1}C presents the truncated 1D Poincar\'e return map $T:~z_0 \to \bar z_{0}$
\begin{equation} 
\bar z_{0}=\mu-B_0R^{1-\nu_0} z_0^{\nu_0}\sin(\Omega_0\ln{z_0}+\phi_2)
\label{1dmap}
\end{equation}
 at $\mu=0$  and with $z>0$. The shape of the map is due to the $\sin$-wave function, with its amplitude or envelope bounded by $\pm B_0R^{1-\nu_0} z_0^{\nu_0}$ ($0<\nu_0<1$), while the frequency of its
 zeros increases logarithmically as $z$ approaches $0^+$. The return map for $z<0$ is the mirror reflection of the above map~\ref{1dmap} . Small variations of $\mu$ vertically shift the map's graph slightly (Fig.~\ref{fig7},~\ref{fig8}), while $z_0$-variations squeeze or stretch it horizontally. 

\begin{figure*}[t!]
\begin{center}
\includegraphics[width=.85\linewidth]{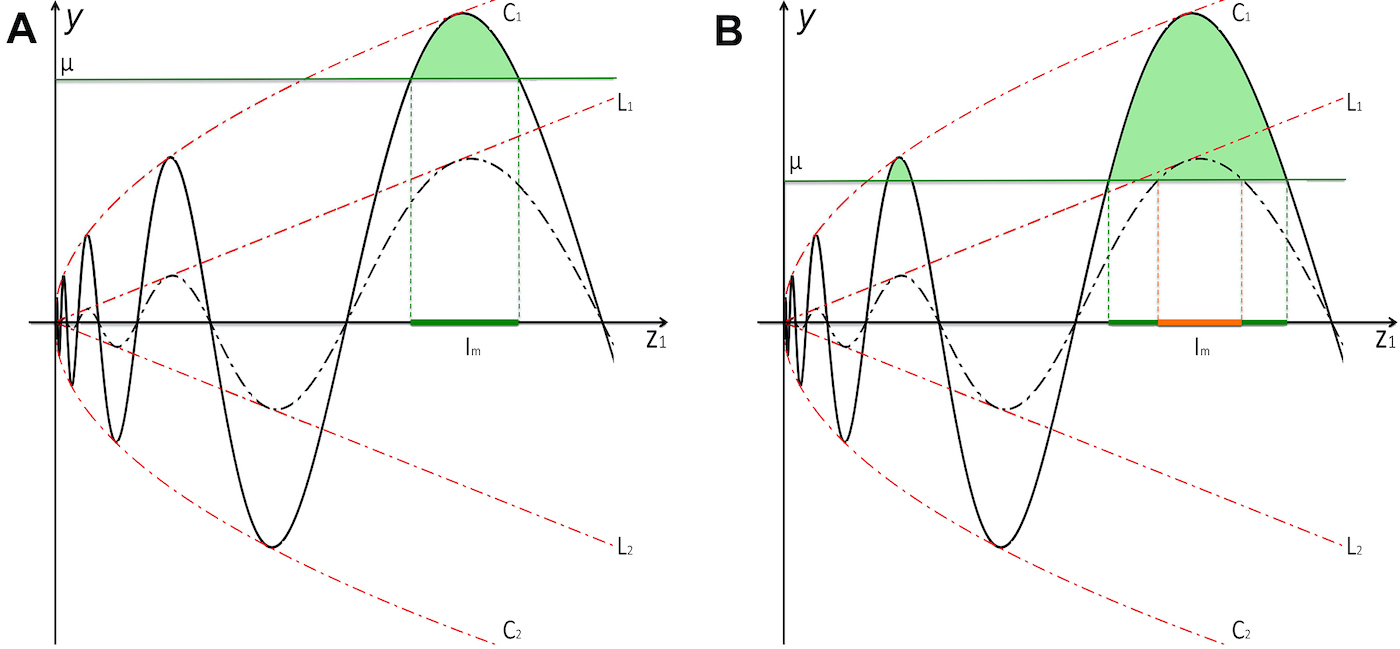}
\caption{The graph (solid line) of the sine-function  $y=B_0R^{1-\nu_0} z_1^{\nu_0}\sin(\Omega_0\ln z_1+\phi_2)$ bounded by the $C_{1,2}$-curves/envelopes. Its section above the horizontal line $y=\mu$ is painted green. The dashed graph is given by the sine-function $y=B_0z_1^{\nu_0}\sin(\Omega_0\ln z_1+\phi_2)$ bounded by the $L_{1,2}$-lines. The dashed graph is completely below the line $y=\mu$ in (A), whereas its crossing with the line $y=\mu$ occurs on an orange interval nested within the green interval $I_m$ on the $z_1$-axis in (B).}\label{fig9}
\end{center}
\end{figure*}

Figure~\ref{fig6} demonstrates several 1D return maps for different $\mu$-values. In particular,  Figs.~\ref{fig6}A and B$_{1,2} $ illustrate multiple [11]- and [00]-homoclinic orbits for a fixed value of the saddle index $\nu_0$, while  Figs.~\ref{fig6}C represents terminal triple [111] and [000] orbits (compare with  Fig.~\ref{fig4}), In these maps, forward iterates of the origin take it back (to $O$) after 2 or 3 steps, respectively. One can infer that such orbits should come in pairs and that the corresponding bifurcation  curves must be of a U-shape, with a critical point associated with tangency, like one in Fig.~\ref{fig6}C. Moreover, an examination of the maps suggests that there are countably many such orbits and bifurcations accumulating to $\mu=0$. Similarly, Figs.~\ref{fig6}D and E show the 1D return maps corresponding to symmetric homoclinic orbits such as [01]/[10] and [110]/[001] that also come in pairs on every oscillation of the map. These conjectures will be proved analytically below.

\subsection{ \bf{[11]}-double homoclinic orbits and bifurcations}

Figure~\ref{fig4} illustrates a [11]-double homoclinic orbit where the 1D unstable separatrix $\Gamma_1$ slightly misses the saddle-focus above its 2D manifold $W^s_{loc}$ after the first loop.  It comes back by intersecting the cross-section $\Pi_0$ at $(\varphi_0,z_0)=(a_1\mu,\,\mu)$, where $\mu>0$. According to the map (\ref{ch7}), the next intersection point of $\Gamma_1$ with $\Pi_0$ is found from these equations
\begin{equation}  
\begin{array}{lcl}
\varphi_1 &=& a_1\mu+A_0R^{1-\nu_0} \mu^{\nu_0}\cos(\Omega_0\ln\mu+\phi_1)+O(\mu^{2\nu_0}), \\
z_1 &=& \mu-B_0R^{1-\nu_0} \mu^{\nu_0}\sin(\Omega_0\ln\mu+\phi_2)+O(\mu^{2\nu_0}),
\end{array} 
\label{ch9}
\end{equation}
where the small term $\varphi_0=a_1\mu$ can be omitted. 

The [11]-double homoclinic orbit occurs when $z_1=0$, i.e.,
\begin{equation}
\mu+O(\mu^{2\nu_0})=B_0R^{1-\nu_0} \mu^{\nu_0}\sin(\Omega_0\ln\mu+\phi_2),
\label{ch10}
\end{equation}
which is equivalent to 
\begin{equation}
\mu^{1-\nu_0}+O(\mu^{\nu_0})=B_0R^{1-\nu_0}\sin(\Omega_0\ln\mu+\phi_2). 
\label{chu1}
\end{equation}

As long as  $\mu$ is  sufficiently small and $\nu_0<1$, then we can assume $B_0R^{1-\nu_0}\sin \left (\Omega_0\ln(\mu)+\phi_2\right )=0$. The solutions of this equation are $\mu_1^{(n)}=e^{-\frac{2n\pi}{\Omega_0}-\frac{\phi_2}{\Omega_0}}$ and $\mu_2^{(n)}=e^{-\frac{2n\pi}{\Omega_0} \pm \frac{\pi}{\Omega_0}-\frac{\phi_2}{\Omega_0}}$ for sufficiently large $n$. In this expression, '+' is used for $B_0>0$, and '-' for $B_0<0$. Without loss of generality, $B_0>0$ is assumed. Note that if $z_1<0$ for $\mu \in (\mu_1^{(n)}, \mu_2^{(n)})$, then $\Gamma_1$ goes underneath $W^s_{loc}$ of the saddle-focus after the second loop, and therefore no sequential one-sided [11....] homoclinic orbits or bifurcations can occur when $\mu \in (\mu_1^{(n)}, \mu_2^{(n)})$.  This situation is illustrated by the 1D maps presented  in Figs.~\ref{fig6}B$_{1,2}$, between which the separatrix falls down below the two zeros of the given U-shaped section of the oscillatory return map.

\begin{figure*}[t!]
\begin{center}
\includegraphics[width=.9\linewidth]{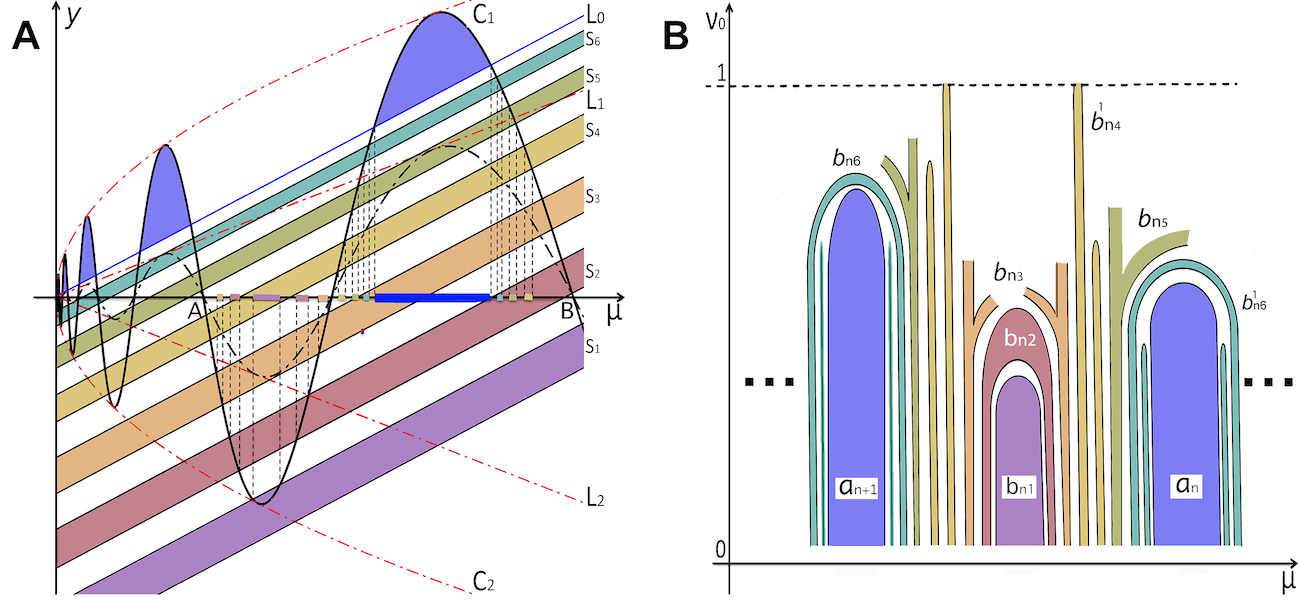}
\caption{The case $B_0<1$. (A) The solid graph of the sine-function $y=B_0R^{1-\nu_0}\mu^{\nu_0}\sin(\Omega_0\ln\mu+\phi_2) $ bounded by the curves $C_{1,2}: y=\pm B_0R^{1-\nu_0}\mu^{\nu_0}$. The dashed graph of the sine function $y=B_0\mu\sin(\Omega_0\ln\mu+\phi_2)$ bounded by the lines $L_{1,2}: y=\pm B_0 \mu$. The interceptions of the solid color strips $S_1-S_6$, parallel to the line $L_0: y=\mu$, and the sine-function are projected to the $\mu$-axis. The closer a strip is to $L_0$, the narrower it becomes. (B) $(\mu,\,\nu_0)$-parameter plane sketching a bifurcation unfolding including  $\cap$-bars (color matching Panel~A) such as $a_n$ and $a_{n+1}$ (in blue) outlining the double [11]/[00] homoclinic orbits, and $b_{nj}$ for triple [111]-orbits.  Around each blue bar there are a pair of narrow $\cap$-bars and a bridge due to the strip $S_6$ in (A). Besides, shown are a few Y-shaped bifurcation objects due to intersections of the the sine-function with the strips $S_3$ and $S_5$, and a pair of yellow $\cap$-bars and a pair of $\Pi$-bars generated by strip $S_4$, as well as a reddish bridge due to the crossing with $S_2$. Shown in the middle is a purple $\cap$-bar $b_{n1}$, which is the widest domain for triple homoclinics orbits; other such regions fit narrow in the parameter space.}\label{fig10}
\end{center}
\end{figure*}

Figures~\ref{fig7}B  and \ref{fig8}B depict the organization of [11]-homoclinic bifurcation curves in the $(\mu,\nu_0)$-parameter plane. Vertical $\cap$-shaped (rounded) bars filled with blue color represent the parameter regions subject to the condition $z_1<0$, whereas their borderlines correspond to the [11]-homoclinic bifurcations, i.e., they are the corresponding bifurcation curves of the [11]-orbits.   The widths $w_n$ and the distances $d_n$ between any two closest $\cap$-shaped bars, evaluated as  $w_n=\mu_2^{(n)}-\mu_1^{(n)}$ and $d_n=\mu_1^{(n)}-\mu_2^{(n-1)}$, resp., decrease proportionately as $n \rightarrow \infty$ ($\mu \rightarrow 0$), as given by the following ratios: 
$$\frac{w_{n+1}}{w_n}=\frac{\mu_2^{(n+1)}-\mu_1^{(n+1)}}{\mu_2^{(n)}-\mu_1^{(n)}}=e^{-\frac{2\pi}{\Omega_0}}
$$ and
 $$\frac{d_{n+1}}{d_n}=\frac{\mu_1^{(n+1)}-\mu_2^{(n)}}{\mu_1^{(n)}-\mu_2^{(n-1)}}=e^{-\frac{2\pi}{\Omega_0}}.
 $$ 
 
Hence, both distance and width shrink exponentially fast as they accumulate to the primary homoclinic bifurcation -- the vertical line $\mu=0$ in the $(\mu,\nu_0)$-diagram (see Figs.~\ref{fig7}B  and \ref{fig8}B).

When $\nu_0 \to 1$, the term $\mu^{1-\nu_0}$ is no longer negligible but significant in Eq.~(\ref{chu1}). In this case, solving $\mu^{1-\nu_0}=B_0R^{1-\nu_0}\sin(\Omega_0\ln\mu+\phi_2)$ or equivalently $\mu=B_0R^{1-\nu_0}\mu^{\nu_0}\sin (\Omega_0\ln\mu+\phi_2)$ gives $z_1=0$. These equations can also be interpreted geometrically, see Figs.~\ref{fig7}A and \ref{fig8}A, respectively, for the cases $B_0<1$ and $B_0>1$.  Namely, the sought condition $z_1=0$  needed for a [11]-homoclinic orbit to close is fulfilled at all intersections of the ``sine'' function $y=B_0R^{1-\nu_0} \mu^{\nu_0}\sin(\Omega_0\ln\mu+\phi_2)$ and the line $L_0:~y=\mu$. As such, the union of all blue intervals (Figs.~\ref{fig7}A and \ref{fig8}A) gives the range of $\mu$-values for which $z_1<0$ for a given constant $\nu_0$.  As $\nu_0 \rightarrow 1$, the amplitudes/envelopes $C_{1,2}: y =\pm B_0R^{1-\nu_0} \mu^{\nu_0}$ of the sine-function $B_0R^{1-\nu_0}\mu^{\nu_0}\sin (\Omega_0\ln\mu+\phi_2)$  flatten and transform into the two lines $L_{1,2}$ given by $y=\pm B_0\mu$ eventually.

If $B_0<1$, the blue $\mu$-intervals start shrinking and vanish after the local maximums of the sine-function are lowered below the line $L_0$, see Fig.~\ref{fig8}A. The closer such a $\mu$-interval is placed to $\mu=0$, the larger value of $\nu_0$ is needed for the interval to vanish. In the $(\mu,\nu_0)$-parameter diagram, the corresponding region looks like a vertical $\sqcap$-bar with the tipping point cut out when $ \mu=B_0R^{1-\nu_0} \mu^{\nu_0}$ or $\nu_0=1-\frac{\ln{B_0}}{\ln{\mu}-\ln{R}}$. This equality is  held on the red dash $C_3$-curve in the $(\mu,\nu_0)$-parameter diagram in Fig.~\ref{fig7}B. This (cusp-shaped) $C_3$-curve approaches the level $\nu_0=1$ from below  as $\mu \rightarrow 0$. The vertical bars that terminate before reaching the horizontal line $\nu_0=1$ all have the $\cap$-shape. 

In the case $B_0>1$, the bars and the intervals between them become narrower as  $\nu_0 \to 1$ but they persist (Fig.~\ref{fig8}A). In the $(\mu,\nu_0)$-parameter diagram, the condition $z_1<0$ is held in the union of all blue $\sqcap$-shaped bars below the level $\nu_0=1$, see Fig.~\ref{fig8}B. 

\subsection{ \bf{[10]}-double homoclinic bifurcations}

\begin{figure*}[t!]
\begin{center}
\includegraphics[width=.85\linewidth]{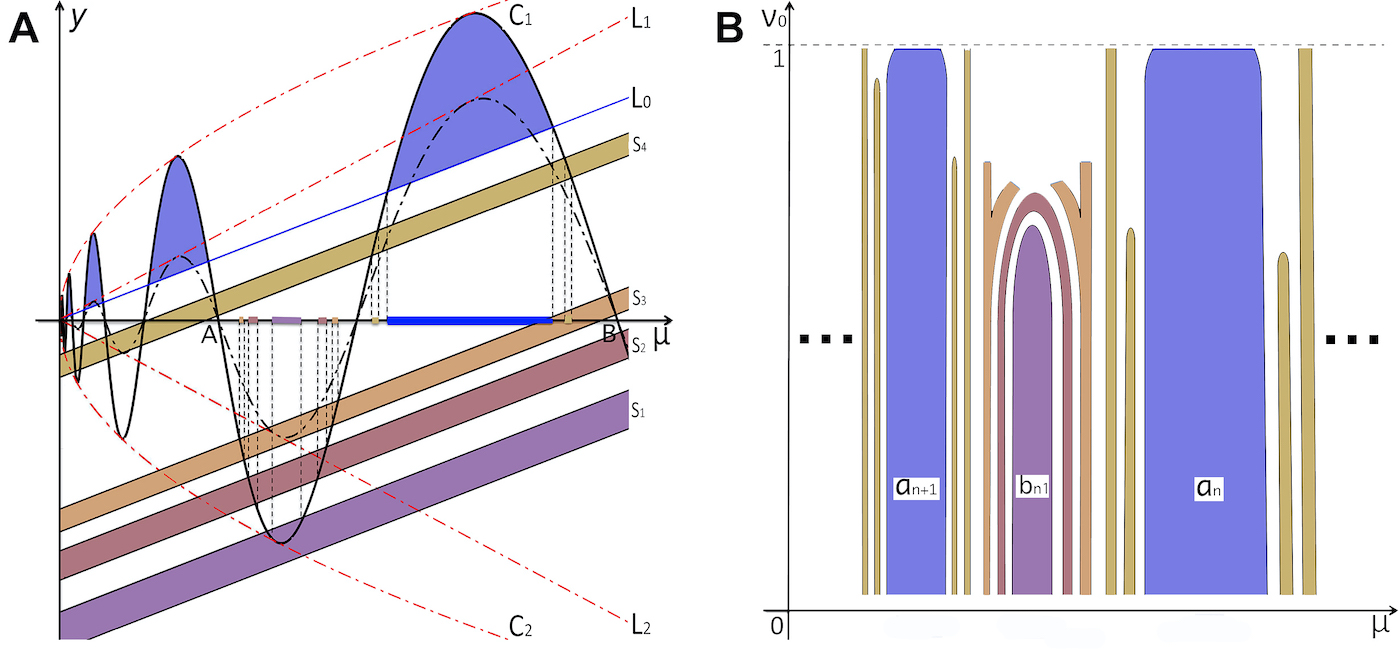}
\caption{$B_0>1$. (a) The sine function is $y=B_0R^{1-\nu_0}\mu^{\nu_0}\sin(\Omega_0\ln\mu+\phi_2) $ bounded by $C_{1,2}: y=\pm B_0R^{1-\nu_0}\mu^{\nu_0}$. The dashed sine function is bounded by $L_{1,2}: y=\pm B_0\mu$. The interceptions of the strips $S_1-S_4$ and the sine function is projected to the $\mu$-axis. The strips are parallel to the line $L_0: y=\mu$. The closer a strip is to $L_0$, the narrower it is. (b) Two blue $(1,1)$ double-loop bars, $a_n$ and $a_{n+1}$ are drawn on the $(\mu,\nu_0)$-plane. The $(1,1,1)$ triple-loop bars fit in between and match the color of the strips in the plot (a) that they are derived from. We describe it in the order from the close region of a blue bar to its furthest. Close to the blue bar, there are a pair of yellow $\cap$-bars and a pair of yellow $\Pi$-bars that are derived from the strip $S_4$. Then there is a pair of "Y" shapes that are derived from $S_3$. Then there is a reddish bridge that is derived from $S_2$. In the middle of the two blue bars, there is a purple $\cap$-bar $b_{n1}$ that is the widest triple-loop piece. The closer to the blue bar, the narrower the triple-loop piece is. Taking $b_{n1}$ as the center, the left pieces are slightly narrower than the right pieces.  }\label{fig11}
\end{center}
\end{figure*}

A typical [10]-double homoclinic orbit is illustrated in Fig.~\ref{fig4}C. The corresponding 1D return map is shown in Fig.~\ref{fig6}D. By construction, after the separatrix $\Gamma_1$ runs a single [1]-loop on its way back to the saddle-focus, it goes underneath its stable manifold $W^s_{loc}$ and hits the cross-section  $\Pi_0$ at $(\varphi_0,z_0)=(a_1\mu,\mu)$ with $\mu<0$, and then it  completes the second [0]-loop heading toward the equilibrium state. It hits $\Pi_0$ for the second time at some point $(\varphi_1, z_1)$, which can be found by the return map as follows:
\begin{equation}  
\begin{array}{lcl}
\varphi_1 &= &\pi+a_1\mu-A_0 R^{1-\nu_0} (-\mu)^{\nu_0}\cos(\Omega_0\ln(-\mu)+\phi_1)+\\&&+O((-\mu)^{2\nu_0}),\\
z_1 &= &-\mu-B_0 R^{1-\nu_0} (-\mu)^{\nu_0}\sin(\Omega_0\ln(-\mu)+\phi_2)+\\&&+O((-\mu)^{2\nu_0}).
\end{array} 
\label{ch11}
\end{equation}
The condition $z_1=0$, i.e., 
$$-\mu+O((-\mu)^{2\nu_0})=B_0 R^{1-\nu_0} (-\mu)^{\nu_0} \sin (\Omega_0 \ln (-\mu)+\phi_2),$$ 
corresponds to the occurrence of a [10]-double homoclinic orbit. One can observe that this condition is similar to the case of [11]-homoclinic orbits where $\mu$ is replaced with $-\mu$. Therefore, the structure of the bifurcation unfolding for [10]-homoclinic orbits is flip-symmetric ($\mu \to -\mu$) to the bifurcation diagram for the [11]/[00]-homoclinics, see Figs.~\ref{fig7}B and \ref{fig8}B.

Now we arrive at the following theorem for double homoclinic orbits. 
\begin{theorem}
Let a reflection-symmetric system have a pair of primary homoclinics to the Shilnikov ($\nu(\mu)<1$) saddle-focus at $\mu=0$.  Then, double homoclinics occur at values   
$$|\mu|=B_0 R^{1-\nu_0} |\mu|^{\nu_0}\sin(\Omega_0\ln |\mu|+\phi_2),$$
where $B_0,\, R>0$ are constants. If $B_0<1$, then the $(\mu,\nu_0)$-parameter diagram includes countably many  $\cap$-shaped bifurcation curves corresponding to double homoclinics that are topped up by the curve, $\nu_0=1-\frac{\ln{B_0}}{\ln{|\mu|}-\ln{R}}$, converging to $1$ as $\mu \to 0$. All such double-bifurcations accumulate to the primary one from both sides with a scalability ratio $e^{-\frac{2\pi}{\Omega_0}}$ for both the width and the distance between the bifurcation curves.
\end{theorem}

\subsection{\bf{[111]}-triple homoclinic orbits and  bifurcations}

A typical [111]-triple homoclinic orbit of the saddle-focus in the phase space is illustrated in Fig.~\ref{fig4}B. Recall that such orbits, [111..], are all one-sided, so to say. The corresponding 1D return map for a critical [111]-orbit is shown in Fig.~\ref{fig6}C where the origin is taken back to zero, here single or critical,  meaning that such a homoclinic orbit may no longer occur for the given map, should the splitting parameter $\mu$ be increased. This tangency at zero in the map corresponds to the turning point of a $\cap$-shaped homoclinic bifurcation curve, like ones shown in Fig.~\ref{fig7}B.   

\begin{figure*}[t!]
\begin{center}
\includegraphics[width=.85\linewidth]{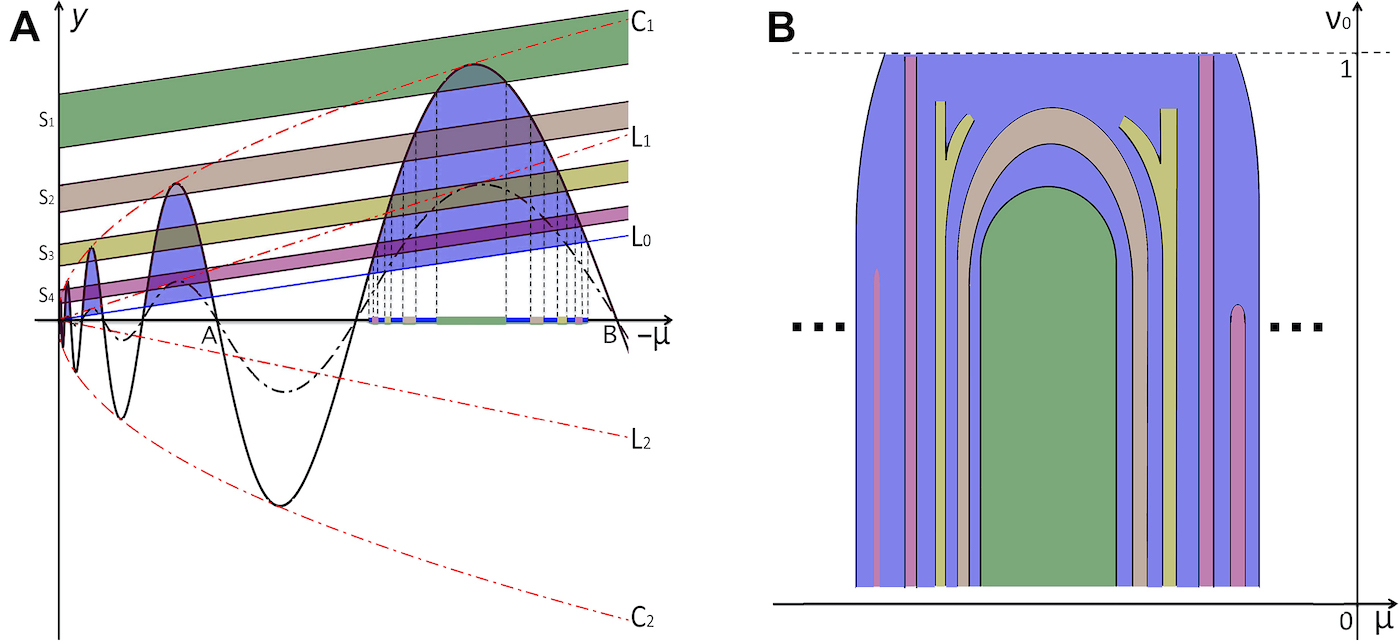}
\caption{$B_0>1$. (a) The x-axis is $(-\mu)$-aixs. The solid sine function is bounded by two curves $C_{1,2}$ and the dashed sine function is bounded by two lines $L_{1,2}$. The parts below solid sine function and above the line $L_0$ are painted blue. Four strips $S_1-S_4$ line up above $L_0$. Their interceptions with the solid sine function are projected to the $(-\mu)$-axis. (b) Inside the blue $Pi$-bar, $(1,0,0)$ triple-loop pieces are drawn to match the strips in the plot (a) that they are derived from. From the middle to both the sides of the blue bar, it lists a green bar, a brown bridge, a pair of yellow "Y" pieces, a pink $Pi$-bar and a pink $\cap$-bar. The width of the green bar in the middle is much larger than the rest. The width of the pieces is decreasing from the middle green bar to the side boundaries of the blue bar. The pieces on the right of the green bar are slightly wider than the pieces on the left.}\label{fig12}
\end{center}
\end{figure*}

So, let the unstable separatrix $\Gamma_1$ make three such loops prior to its returning to the saddle-focus along its $W^s_{loc}$. Then, the corresponding truncated map, accounting for the dominant terms only, is given by 
\begin{equation}  
\begin{array}{lcl}
z_1&=&\mu-B_0R^{1-\nu_0} \mu^{\nu_0}\sin(\Omega_0\ln\mu+\phi_2)+O(\mu^{2\nu_0}),\\
z_2&=&\mu-B_0R^{1-\nu_0} z_1^{\nu_0}\sin(\Omega_0\ln z_1+\phi_2)+O(z_1^{2\nu_0})\\
\end{array} 
\label{ch12}
\end{equation}
(equations for $\varphi$-variables  are omitted). To find the bifurcation curves corresponding to  [111]-triple homoclinics, one must first identify the $\mu$-range where $z_2<0$.

The first part in Eqs.~(\ref{ch12}) implies that $z_1\sim\mu^{\nu_0}$, and therefore $z_1$ is small when $\mu$ is small. The second equation in (\ref{ch12}) can be further reduced to $\sin(\Omega_0 \ln z_1+\phi_2)=0$, assuming that $\mu$ is small enough. Its solutions are $z_1=e^{-\frac{2m\pi}{\Omega_0}-\frac{\phi_2}{\Omega_0}}$ or $e^{-\frac{2m\pi}{\Omega_0} + \frac{\pi}{\Omega_0}-\frac{\phi_2}{\Omega_0}}$, here $m \in \mathbb{Z} $ is sufficiently large. Denote by $I_m$ all intervals of $z_1$ for such $z_2<0$. It can be deduced from the $m$-th period of the sine-function that the sought intervals are given by 
$$
I_m \approx (e^{-\frac{2m\pi}{\Omega_0}-\frac{\phi_2}{\Omega_0}}, e^{-\frac{2m\pi}{\Omega_0} + \frac{\pi}{\Omega_0}-\frac{\phi_2}{\Omega_0}})
$$ 
with $m \in \mathbb{N}$ being sufficiently large, and $\mu$ is assumed to be small. Figure~\ref{fig9}  illustrates such an interval $I_m$ for some small fixed $\mu$: it is highlighted in green on the $z_1$-axis within which the $m$-th period of the sine-function, $y=B_0R^{1-\nu_0} z_1^{\nu_0}\sin(\Omega_0\ln z_1+\phi_2)$, is greater than the given  $\mu$. Assume that $D_m$ is a local maximum on the $m$-th period of the sine-function whose graph is a dashed-line in 
Fig.~\ref{fig9}B. If $D_m<\mu$, then any such an interval $I_m$ shrinks and collapses as $\nu_0 \rightarrow 1$. If $D_m>\mu$, in contrast, then $I_m$ narrow downs to some (orange) interval that persists in the limit $\nu_0 \rightarrow 1$, see Fig.~\ref{fig9}B.

\begin{figure*}[t!]
\begin{center}
\includegraphics[width=.85\linewidth]{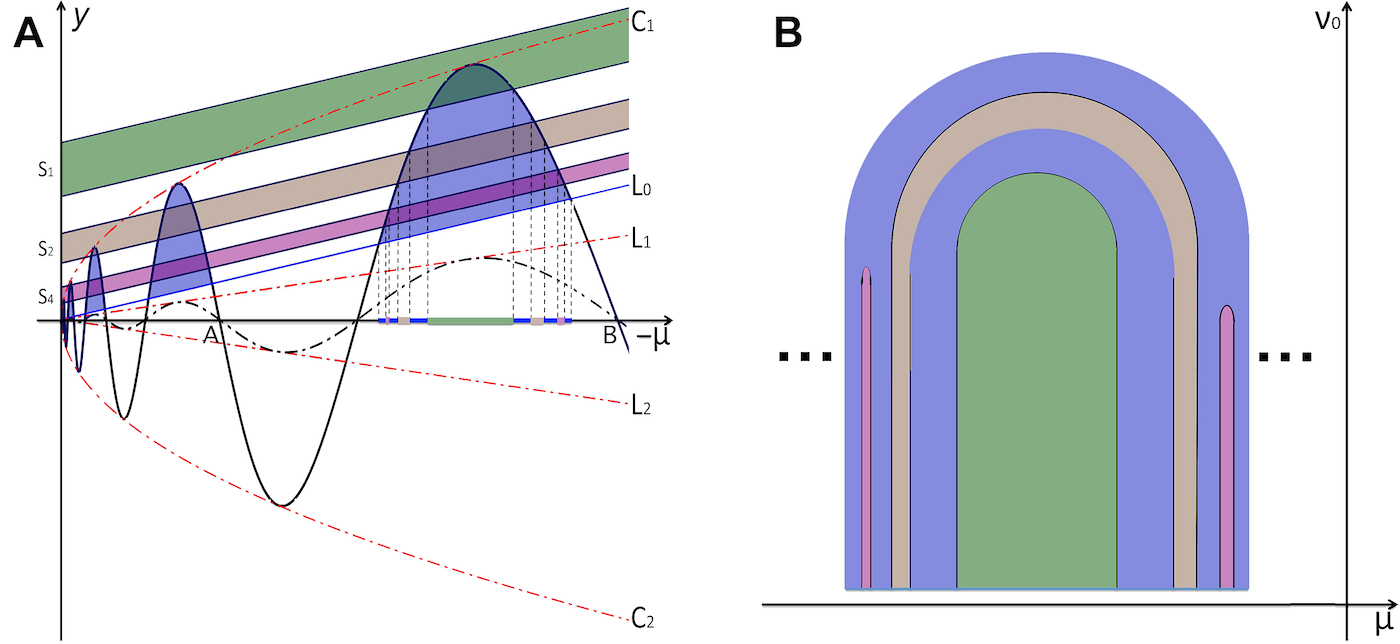}
\caption{$B_0<1$. (a) The x-axis is $(-\mu)$-axis. The solid sine function is bounded by two curves $C_{1,2}$ and the dashed sine function is bounded by two lines $L_{1,2}$. The parts below solid sine function and above the line $L_0$ are painted blue. Three strips $S_1$, $S_2$ and $S_4$ line up above $L_0$. Their interceptions with the solid sine function are projected to the $(-\mu)$-axis. (b) Inside the blue $\Pi$-bar, $(1,0,0)$ triple-loop pieces are drawn to match the strips in the plot (a) that they are derived from. From the middle to both the sides of the blue bar, it lists a green bar, a brown bridge, and a pink $\cap$-bar. The width of the green bar in the middle is much larger than the rest. The width of the pieces is decreasing from the middle green bar to the side boundaries of the blue bar. The pieces on the right of the green bar are slightly wider than the pieces on the left. }\label{fig13}
\end{center}
\end{figure*}

Let us first discuss the case $B_0<1$. The first equation (\ref{ch12}) can be written as follows
$$
\mu-z_1=B_0R^{1-\nu_0} z_1^{\nu_0}\sin(\Omega_0\ln z_1+\phi_2)+O(z_1^{2\nu_0}).
$$
One can observe that the right hand-side of the equation above is the sine-function, like one shown in Fig.~\ref{fig10}A. Its sections above $L_0$, which are filled in blue, are the same ones shown in Fig~\ref{fig7}A, and the boundaries of the blue $\cap$-shape bars in Fig.~\ref{fig10}B are the [11]-homoclinic bifurcation curves (Fig~\ref{fig7}B), which are elaborated on in the previous section. Let us examine [111]-homoclinic bifurcations occurring only on one period of the sine-function $y=B_0R^{1-\nu_0}\mu^{\nu_0} \sin(\Omega_0 \ln \mu+\phi_2)$. The chosen period is labeled by endpoints, A and B, in Fig.~\ref{fig10}A. Consider the $n$-th period given by $\mu \in (e^{-\frac{2n\pi}{\Omega_0}-\frac{\pi}{\Omega_0}-\frac{\phi_2}{\Omega_0}}, e^{-\frac{2n\pi}{\Omega_0} + \frac{\pi}{\Omega_0}-\frac{\phi_2}{\Omega_0}})$. Assume that, the function has a local minimum $m_n$ at $\mu_{min}$, and a local maximum $M_n$ at $\mu_{max}$. After  the bent envelops $C_{1,2}$ are rectified and become the straight lines $L_{1,2}$ at $\nu_0=1$, the new local minimum and maximum on the given period are denoted by $m_n'$ and $M_n'$. 
%Without designation, the sine-function we mention below has this n-th period only. 
Consider [111]-homoclinic bifurcations occurring between two blue bars, $a_n$ and $a_{n+1}$, (corresponding to the occurrence of two consecutive [11]-orbits on the same  or similar interval) in the bifurcation diagram in Fig.~\ref{fig10}B.

Recall that $z_2<0$ when $z_1 \in I_m$, with $m \in \mathbb{N}$ being sufficiently large. 
By construction,  $y=\mu-I_m$ represents a strip bounded by two lines nearly parallel  for large enough $m \in \mathbb{N} $, see Fig.~\ref{fig10}; actually, any two successive lines  are no longer parallel if we take into account smaller terms that were neglected earlier.   Therefore, the range of $\mu$-values corresponding to $z_2<0$ is a union of all the interception intervals of (colored) strips $y=\mu-I_m$,  with the curve $y=B_0R^{1-\nu_0}\mu^{\nu_0} \sin(\Omega_0 \ln \mu+\phi_2)$ projected on the $\mu$-axis,  All such strips line up under and accumulate from below to $L_0$ as $m$ increases. They become narrower while  approaching $L_0$ so that $\frac{|I_{m+1}|}{|I_m|}=e^{-\frac{2\pi}{\Omega_0}}$.  Six such colored strips labeled by $S_1\cdots S_6$ are sampled in Fig.~\ref{fig10}A corresponding to multiple distinct [111]-triple orbits and their homoclinic bifurcation curves in the $(\mu,\nu_0)$-parameter plane (Fig.~\ref{fig10}B). 

Let the purple strip $S_1$ at the bottom be the very first one that intercepts the oscillatory graph of the sine-function $y=B_0R^{1-\nu_0}\mu^{\nu_0} \sin(\Omega_0 \ln \mu+\phi_2)$. The projection of this overlap onto the $\mu$-axis, say $p_1$, is the $\mu$-interval where $z_2<0$ for a given constant $\nu_0$. It is located in between two blue intervals, on which the sine function is greater than $L_0$ (Fig.~\ref{fig10}A). In the limit $\nu_0 \rightarrow 1$, the bended envelopes $C_{1,2}$ straighten up and become the $L_{1,2}$-lines, which makes the given $p_1$-interval collapse and vanish. The corresponding image of $p_1$  is the purple $\cap$-shaped solid bar, say $b_{n1}$, located in the middle of the two blue $\cap$-shaped bars, $a_n$ and $a_{n+1}$ (corresponding to the condition $z_1<0$) in the $(\mu,\,\nu_0)$-parameter plane in Fig.~\ref{fig10}B. In the case of $S_2$ or $S_6$, we have the following inequalities $m_n<(\mu-I_m)|_{\mu=\mu_{min}}<m_n'$ or $M_n'<(\mu-I_m)|_{\mu=\mu_{max}}$, respectively, i.e., the stripes are bounded by the old and new  local minima and maxima on the $n$-th period of the sine-function. Fig.~\ref{fig10}B  gives an interpretation of these inequalities in the $(\mu,\,\nu_0)$-bifurcation diagram: the corresponding (reddish) region, say $b_{n2}$, is formed through a merger of two bending vertical bars forming a bridge- or arch-like connection atop of the $\cap$-shaped one $b_{n1}$; same is true for the green bending bridges, say $b_{n2}$ and $b_{n6}^1$ due to $S_6$-strip, which are placed on top of the blue bars, $a_n$ and $a_{n+1}$. The geometric explanation of such $\cap$-shape is the same: graph of the sine-function no longer crosses $S_2$ or $S_6$ as $\nu_0$ is increased or decreased beyond some thresholds corresponding to critical tangencies with these stripes. Clearly, there can be more such bridge-shaped regions in the bifurcation diagram  if there are more such strips fitting these conditions.

In addition, a top strip like $S_6$ can give rise to extra narrow (green) $\cap$-bars located inside every  bridge (such as $b_{n6}$ and  $b_{n6}^1$) in Fig.~\ref{fig10}B provided that  it is (i) close to $L_0$ and (ii) is narrow enough to shrink and vanish before the corresponding two bending bars merge, forming a bridge-like object in the bifurcation diagram in the limit $\nu_0 \rightarrow 1$ .  Otherwise, if the following inequalities $m_n'<(\mu-I_m)|_{\mu=\mu_{min}}$ and $(\mu-I_m)|_{\mu=\mu_{max}}<M_n'$ are fulfilled, say for the  $S_4$-strip in Fig.~\ref{fig10}A, then no bridge- or arch-shaped region will be formed through a merger of two bending bars. Instead, the corresponding (yellow) structure will look like either a $\Pi$-bar ($b_{n4}^1$) or a $\cap$-bar ($b_{n4}^2$ to the right from $b_{n4}^1$), in the $(\mu,\,\nu_0)$-parameter plane in Fig.~\ref{fig10}B. The shape of such bars, $\Pi$- or $\cap$-like, is determined by whether the width of the strip or, equivalently, the width of the $I_m$-interval remains small but finite as in Fig.~\ref{fig9}B, or it collapses as depicted in Fig.~\ref{fig5}A in the limit $\nu \rightarrow 1$. Since a narrower strip is likely to vanish, therefore, the yellow $\cap$-bar resides  closer to the blue [11]-orbit bars. If $m_n'$ or $M_n'$ happen to be inside the given strip, like $S_3$ or $S_5$, then two of its sections intercepting the graph of the sine-function will merge after the minimum or maximum of the sine-function move inside $S_3$ and $S_5$, with changes in $\mu$. Meanwhile both $S_3$ and $S_5$ keep narrowing as $\nu_0 \rightarrow 1$. Therefore, the extreme (min/max) points can slip away from overlapping with $S_3$ or $S_5$, which result in the corresponding bridge-like region decoupling into two ``Y''-shaped objects. These are depicted in Fig.~\ref{fig10}B as yellow pair-wise bars, like $b_{n3}$ due to $S_3$, atop of $b_{n2}$, and greenish Y-shaped bars like $b_{n5}$ due to $S_5$.  Both branches remain left-open as they may end up differently with parameter variations. Either branch of a Y-shaped bar can bridge with the symmetric one, or they both terminate prior to merging into one. A bar like $b_{n3}$ can also morph into the shape of the bar next to it, of  the $\Pi$- or $\cap$-shape. Thus, the Y-shape can be viewed as a transition between the bridge and the $\cap$-bar.  It is easy to argue that there can be a single Y-shaped region or none generated by each strip like $S_3$ or $S_5$. Let $b_{n1},\cdots,b_{n6}$ be referred to as [111]-triple homoclinic zones. One or more of such zones would be absent in the bifurcation diagram sketched in Fig.~\ref{fig10}B, if there were no strips passing throughout the corresponding position in  Fig.~\ref{fig10}A and so forth. For example, the $b_{n1}$-zone may no longer be present if the strip $S_1$ were positioned to go through the sine-function in the way the strip $S_2$ does. Another example is $S_4$: if the relative positions of $\mu$ and the $z_1$ sine-function for all $S_4$ strips were such as shown in Fig~\ref{fig9}B, than the $b_{n4}^2$-bars would not be spotlighted in Fig.~\ref{fig10}B.

\begin{figure*}[hbt!]
\begin{center}
\includegraphics[width=.99\linewidth]{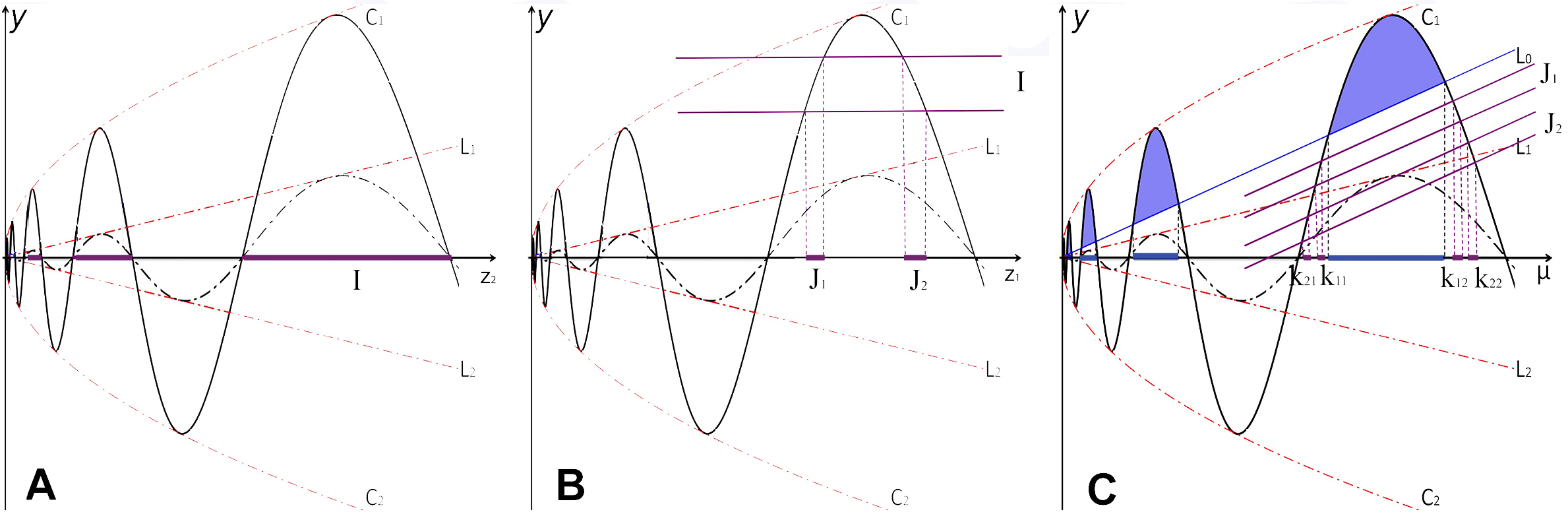}
\caption{(a) The positive parts of the sine function of $z_2$ is projected onto $z_2$-axis and painted purple. One of them is notated as $I$.(b) A purple interval from (a),$I$, corresponds to a strip that intercepts the sine function of $z_1$. The purple intervals on the $z_1$-axis are the projections of the intercepts and they are notated as $J_1$ and $J_2$. (c) The intervals from (b) corresponds to two strips and they intercept the sine function of $\mu$. The purple intervals on the $\mu$-axis are the projections of the intercepts and they are notated as $k_{11}$, $k_{12}$, $k_{21}$ and $k_{22}$ respectively.  }\label{fig14}
\end{center}
\end{figure*}

Note that a [111]-zone in Fig.~\ref{fig10}B becomes the thicker, the further it is away from the closest blue $\cap$-bars $a_n$ and $a_{n+1}$. So, $b_{n1}$, if it exits, is significantly larger than others because the strip $S_1$ intercepts the graph of the sine-function at its flattest section near the critical point. Actually,  by computing the derivative of the sine-function $y=B_0R^{1-\nu_0}\mu^{\nu_0}\sin(\Omega_0\ln\mu+\phi_2)$, we obtain 
\begin{equation}
\frac{d\,y}{d\mu}=
\frac{N_2}{\mu^{1-\nu_0}}\sin(\Omega_0\ln\mu+\phi_2+\theta_0),
\label{ch13}
\end{equation}
where $N_2=B_0R^{1-\nu_0} \sqrt{\nu_0^2+\Omega_0^2}$ and $\cos\theta_0=\nu_0/\sqrt{\nu_0^2+\Omega_0^2}$. For a fixed $\nu_0<1$, $\frac{dy}{d\mu}$ is large because $\mu$ is small. Therefore, the graph of the sine-function looks as if it is made of vertical (and horizontally dense) lines, except for small neighborhoods of critical points of $y=B_0R^{1-\nu_0}\mu^{\nu_0}\sin(\Omega_0\ln\mu+\phi_2)$. The overlapping of $S_1$ with the sine-function is one such neighborhood. Therefore, the width of $b_{n1}$ is to be significantly larger compared to those of other [111]-zones. 

On the other hand, the width of a [111]-zone increases with increasing $\mu$ because $\frac{dy}{d\mu}$ decreases in Eq.~(\ref{ch13}); see such zones on the right from $b_{n1}$ sketched in Fig.~\ref{fig10}B, which should look slightly wider than the symmetric zones on the left from $b_{n1}$.  We conclude by remarking that the very middle zone -- the purple bar $b_{n1}$  being furthest from $a_n$ and $a_{n+1}$ is sketched to be of the largest width. Note that in computational sweeps, such an associated purple bar, or a reddish bridge if the former one does not exist, can be   the only visible or recognizable [111]-zone, as all others might be too narrow to detect. We reiterate that the borderlines of such zones in in the $(\mu,\nu_0)$-parameter plane are the bifurcation curves corresponding to [111]-triple homoclinic orbits.

Figure \ref{fig11} is meant to aid with describing the region where $z_2<0$, and hence with detecting  [111]-homoclinic bifurcations when $B_0>1$. Unlike its predecessor, it does include the strips $S_5$ and $S_6$ and the corresponding bifurcation zones. One can see that the bifurcation diagram in Fig.~\ref{fig11}B, still featuring the bridges atop of $b_n$ and $a_{n+1}$ along with transitional Y-shaped bars, is similar to that in the case $B_0<1$.

\subsection{{\bf [100]}-triple homoclinic orbits and bifurcations}

A similar, [110]-triple homoclinic orbit is pictured in Fig.~\ref{fig4}D.  Unlike it, the unstable separatrix $\Gamma_1$ makes one loop above and two loops underneath the stable manifold $W^s$  before it returns to the saddle-focus to complete a [100]-homoclinic orbit. The corresponding $z$-return map with the dominant terms only is given by 
\begin{equation}  
\begin{array}{lcl}
z_1&=&-\mu-B_0R^{1-\nu_0}(-\mu)^{\nu_0}\sin(\Omega_0\ln(-\mu)+\phi_2)+\\&&+O((-\mu)^{2\nu_0}),\\
z_2&=&-\mu+B_0R^{1-\nu_0}(-z_1) ^{\nu_0} \sin(\Omega_0\ln(-z_1)+\phi_2)+\\&&+O((-z_1)^{2\nu_0}),\\
\end{array} 
\label{ch14}
\end{equation}
where $\mu<0$.  As before, we seek the range of $\mu$-values such that $z_2<0$ for a fixed $\nu_0$. By varying $\nu_0$, we identify the regions  in the $(\mu,\nu_0)$-bifurcation diagram corresponding to the condition $z_2<0$. The boundaries of these regions are the bifurcation curves associated with  [100]-homoclinic orbits.  

\begin{figure*}[t!]
\begin{center}
\includegraphics[width=.8\linewidth]{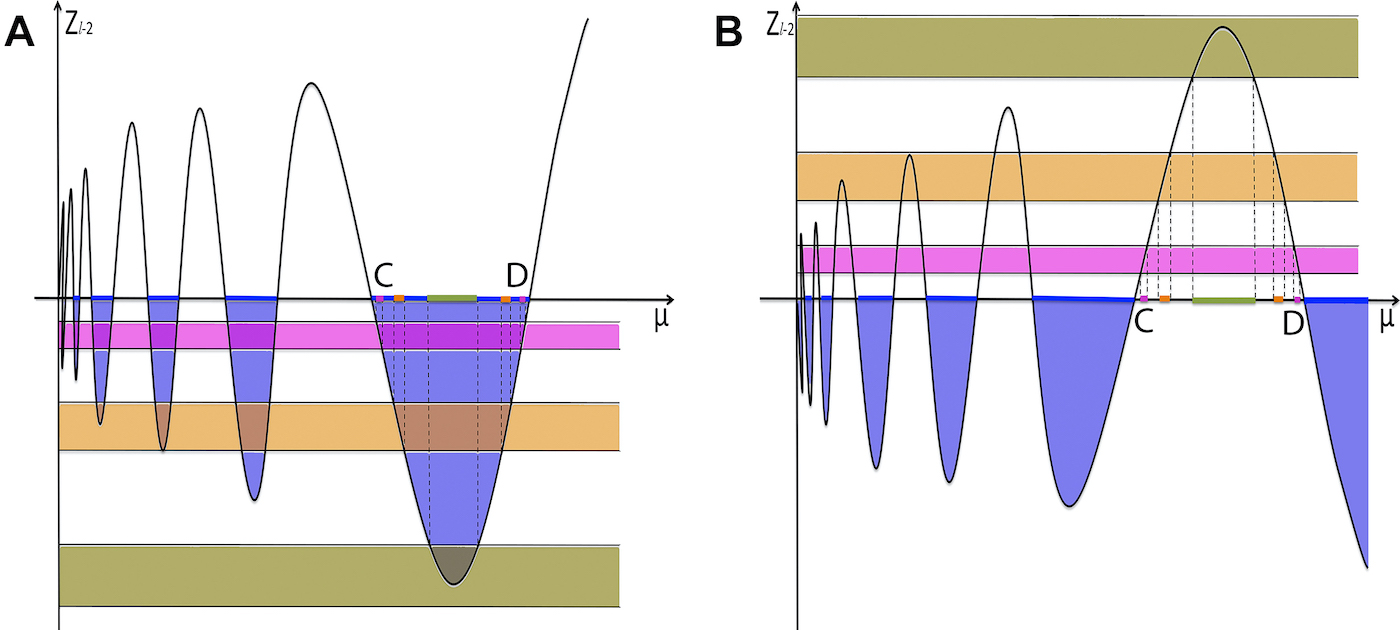}
\caption{The oscillating function is assumed to be $z_{l-2}=F(\mu)$ where $\mu>0$. The regions of $z_{l-2}<0$ are marked blue. Their projections on the $\mu$-axis are marked as blue intervals. Points $C$ and $D$ are the endpoints of one of the intervals.(a) In the case of $\alpha_{l-1}=0$, the strips of $z_{l-1}=f_{\alpha_{l-1}}(z_{l-2})<0$ are beneath the $\mu$-axis. The projected intervals of the interceptions of the strips and the oscillating function are inside the blue intervals. (b) In the case of $\alpha_{l-1}=1$, the strips of $z_{l-1}=f_{\alpha_{l-1}}(z_{l-2})<0$ are above the $\mu$-axis. The projected intervals of the interceptions of the strips and the oscillating function fall in the gaps of the blue intervals.}\label{fig15}
\end{center}
\end{figure*}

  Equations~(\ref{ch14}) imply that  $-z_1\sim(-\mu)^{\nu_0}$, hence $z_1$ is small when $-\mu$ is small. Omitting the small term $-\mu$, the second equation in (\ref{ch14}) can be further simplified: $\sin(\Omega_0\,\ln(-z_1)+\phi_2)=0$, i.e., $-z_1=e^{-\frac{2m\pi}{\Omega_0}-\frac{\phi_2}{\Omega_0}}$ or $e^{-\frac{2m\pi}{\Omega_0} -\frac{\phi_2}{\Omega_0}+\frac{\pi}{\Omega_0}}$ with $m \in \mathbb{Z} $ large enough. The interval $I_m$ of $-z_1$-values such that $z_2<0$, derived from the $m$-th period of the sine-function, can be estimated as $\left (e^{-\frac{2m\pi}{\Omega_0}-\frac{\phi_2}{\Omega_0}},\, e^{-\frac{2m\pi}{\Omega_0} -\frac{\phi_2}{\Omega_0}+\frac{\pi}{\Omega_0}}\right )$. As  was discussed previously, $I_m$ is actually $\mu$-related if we consider some negligible small $\mu$-term. The length of $I_m$ can either decrease to zero (Fig.~\ref{fig9}A), or decrease to a small fixed number (Fig.~\ref{fig9}B) as $\nu_0 \rightarrow 1$. 

\begin{figure*}[t!]
\begin{center}
\includegraphics[width=.9\linewidth]{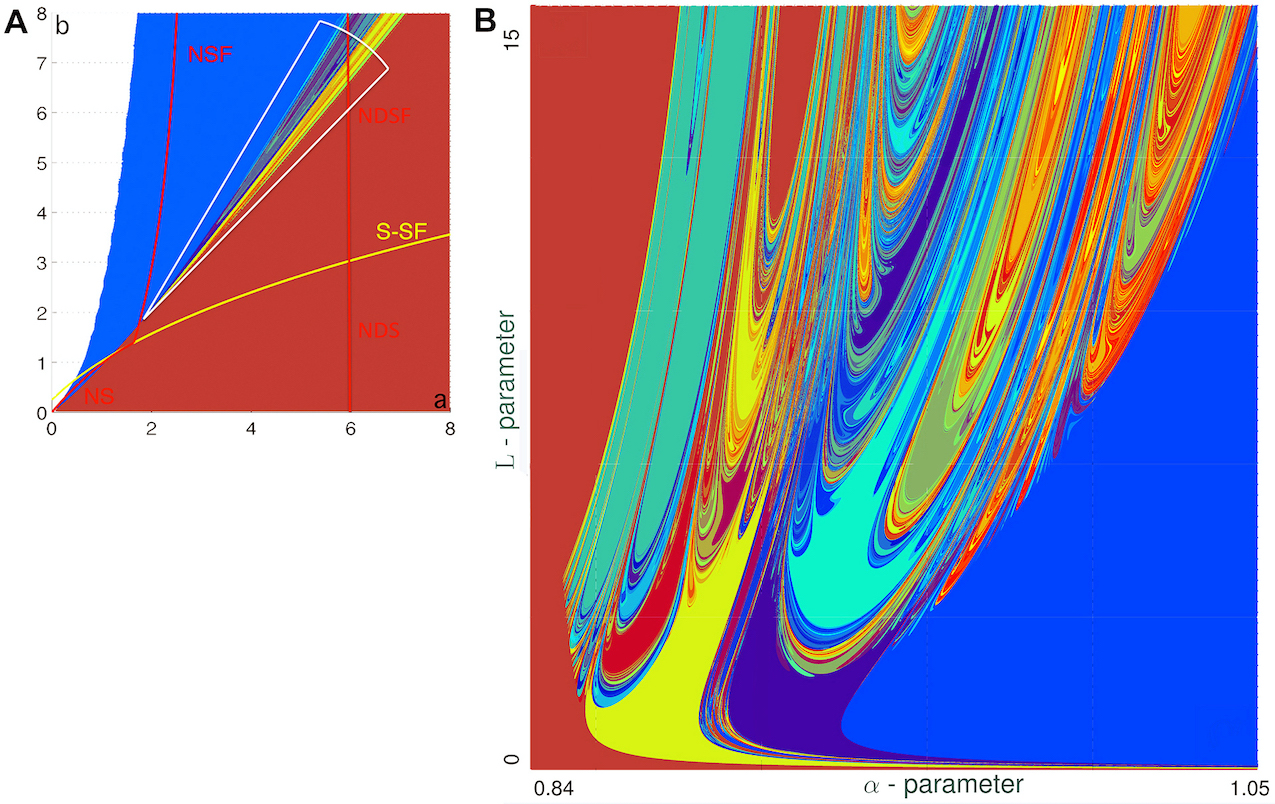}
\caption{(A) Short bi-parameter sweep of the cubic Chua model using the original parameters. The curves NSF and NS correspond to zero saddle value $\sigma_1=0$ or $\nu=1$ ($\nu>/<1$ on the left/right of these curves); NDSF and NDS correspond to zero divergence $\sigma_2=0$ or the saddle index $\nu=1/2$ at the saddle-focus. The yellow S-SF curve marks  the transition of the origin between a saddle and a saddle-focus. All three curves cross near  the point $(1.6458,1.3934)$. The solid color regions correspond to trivial dynamics with constant kneadings. The narrow colorful wedge is the region of double-scroll chaotic dynamics due to  homoclinic bifurcations. The sector bounded by the white curves with the tip point near $(1.8623,1.8743)$ is unfold in Panel~B. (B) Short bi-parameter sweep of the model with new  $(\alpha, L)$-parameters revealing the stunning complexity and universality of homoclinic bifurcations of the Shilnikov saddle-focus at the origin of the Chua model; here, [2-12] and [6-15]-long binary sequences are used for A and B, resp.}\label{fig16}
 \end{center}
\end{figure*}

Let us first elaborate on the case $B_0>1$.  The first equation of the system (\ref{ch14}) can be written as 
$$-\mu-z_1=B_0R^{1-\nu_0}(-\mu)^{\nu_0} \sin(\Omega_0 \ln(-\mu)+\phi_2) +O((-\mu)^{2\nu_0}).$$ 
The right hand-side is the sine-function of $-\mu$, whose graph is depicted in Fig.~\ref{fig12}A, with $-\mu$ being on the x-axis. Its graph sections,  painted in blue, above $L_0:~y=-\mu$ let the $\mu$-intervals be identified within which [10]-homoclinic bifurcations occur, as discussed previously. Recall that $z_2<0$ when $-z_1 \in I_m$ (here, $m \in \mathbb{Z} $ is to be large enough), and therefore the range of $\mu$-values such that $z_2<0$ is the projection of the interceptions of the line $y=(-\mu)+I_m$  with the graph of  $y=B_0R^{1-\nu_0}(-\mu)^{\nu_0}\sin(\Omega_0\ln(-\mu)+\phi_2)$ onto the $(-\mu)$-axis. Then, for each $m$, the quantity $y=(-\mu)+I_m$ is geometrically interpreted as a strip that is parallel to and above the line $L_0$. The domain of $-\mu$ such that $z_2<0$ is the projection of the overlaps of all such strips and the graph of the sine-function onto the $(-\mu)$-axis.  Four such strips $S_1$ -- $S_4$ are  sampled in Fig.~\ref{fig12}A to help us examine [100]-homoclinic bifurcations. Start with the interception of the top (green) strip $S_1$ and the sine-function graph: its projection is an interval on the $(-\mu)$-axis that collapses to zero as $\nu_0\rightarrow 1$ that makes the envelope $C_{1,2}$ converge to the lines $L_{1,2}$. Its image in the  $(\mu,\nu_0)$-diagram shown in Fig. \ref{fig12}B is a $\cap$-shaped (green) bar in the middle of the blue wide $\Pi$-bar corresponding to all homoclinic orbits starting with the [10]-code. The brown strip $S_2$ (located above the the dashed graph of the sine-function) overlaps with the blue zone (on period A--B) on two $\mu$-intervals which merge and then vanish when $\nu_0 \rightarrow 1$. Its image in the bifurcation diagram in Fig~\ref{fig12}B is a (brownish) bridge or arch above the green $\cap$-shaped bar. Let the strip $S_3$ cover the local maximum of the dashed sine-function in Fig~\ref{fig12}B. Such intersection infers respectively that its corresponding images in the bifurcation diagram can be  pair of  Y-shaped branches within  the blue bar, see Fig.~\ref{fig12}B. The pink strip $S_4$ is placed under the local maximum of the dashed sine-function. Therefore, its two overlaps cannot merge, and therefore they correspond to a pair of narrow (pink) $\Pi$-bars or $\cap$-bars in Fig.~\ref{fig12}B. Observe that the $\cap$-shaped bars are located closer to the border of the enclosing blue $\Pi$-bar. As argued previously, there can be only one such green $\cap$-bar in the middle and a single pair of Y-shaped branches, if any, unlike the bars of other shapes for [100]-homoclinics that are not shown in Fig.~\ref{fig12}B to make it visually less busy. Note that the green $\cap$-bar is wider than others because it is due to the overlap of the strip $S_1$ with  the flattest part of the sine-function near its local max. To conclude, let us recap that all [100]-homoclinic bifurcation objects fit inside the [10]-region. As such, one can likely notice regions associated with the green bar for [100]-homoclinics in bi-parametric  sweeps of real applications, as we will demonstrate in the second computational part of our paper. All other bars are probably too slim compared to the principle one, in the given scale, as they originate in the region where the applied sine-function looks as if it is  composed of nearly vertical oscillatory segments. We emphasize that the borderlines of the  bars described and sketched in Fig.~\ref{fig12}B correspond to the [100]-homoclinic bifurcations.

In the case $B_0<1$, the dashed sine-function graph resides fully outside of the blue regions. This case is  somewhat similar to the case of $B_0>1$, except that there are no yellow $\mathbb{Y}$-shaped (due to obvious reasons by their construction) and pink $\Pi$-shaped bars in the bifurcation diagram presented in Fig.~\ref{fig13}B. The pink $\Pi$-like bars (not the $\cap$-shape) cannot exist when $B_0<1$  because the blue $\Pi$-shaped region morphs into the $\cap$-shape. One can deduct from Fig.~\ref{fig13}A that all overlaps and interceptions within the blue domain first shrink, and then disappear after the dashed graph is lowered further below some point.

\subsection{\bf One-sided  [111$\cdots$]--homoclinic orbits and bifurcations}

A one-sided  [111$\cdots$]--homoclinic orbit of the saddle-focus is a longer extension of the  [111]-triple homoclinic orbit as one depicted in Fig.~\ref{fig4}B. The corresponding 1D return maps are somewhat similar to that presented in Fig.~\ref{fig6}C, with the difference that it takes more forward iterates of the origin to come back to zero.      

Let us consider the case where the right unstable separatrix $\Gamma_1$ of the saddle-focus of the origin orbits, say $l$ ($l>3$), one-sided loops before it touches the stable manifold $W^s$. Its symbolic representation is coded as $[111\cdots\,1]_{l}$ . The system of the Poincar\'{e} map accounting for the dominant terms only (equations for $\varphi$s  are omitted) can be written as follows:
\begin{equation}  
\begin{array}{lcl}
z_1&=&\mu-B_0R^{1-\nu_0}\mu^{\nu_0} \sin ( \Omega_0 \ln \mu+\phi_2)+O(\mu^{2\nu_0}),\\
z_2&=&\mu-B_0R^{1-\nu_0} ~~z_1 ^{\nu_0} \sin ( \Omega_0 \ln z_1+\phi_2)+O(z_1^{2\nu_0}),\\
 & \vdots     & \\
z_{l-2}&=&\mu-B_0R^{1-\nu_0} z_{l-3} ^{\nu_0} \sin ( \Omega_0 \ln z_{l-3}+\phi_2)+O(z_{l-3}^{2\nu_0}),\\
z_{l-1}&=&\mu-B_0R^{1-\nu_0} z_{l-2} ^{\nu_0} \sin ( \Omega_0 \ln z_{l-2}+\phi_2)+O(z_{l-2}^{2\nu_0}).\\
\end{array} 
\label{ch15}
\end{equation}
Our goal here is to determine the structure of the corresponding bifurcation curves,  assuming that we have already known all unfoldings for the shorter loops up to order $l-1$. The equivalent problem is to find $R_l=\{\mu|z_{l-1}<0\}$, given that $R_i=\{\mu|z_{i-1}<0\}$ ($2 \le i \le l-1$). Obviously,  $z_{i-1}>0$ ($2 \le i \le l-1$) for $\mu \in R_l$, see Eqs.~(\ref{ch15}) above. Therefore, $R_l$ disjoints $\bigcup_{i=2}^{l-1} R_i$; i.e., $R_l \subset \overline{\bigcup_{i=2}^{l-1} R_i}$. This guarantees that bifurcation curves for the right $l$-loops  fit into the gaps between all other one-sided (left/right) homoclinic orbits of lower orders. 

System (\ref{ch15}) can also be recast as
 $$
 \begin{array}{lcl}
 z_1&=&\mu-B_0R^{1-\nu_0}\mu^{\nu_0} \sin ( \Omega_0 \ln \mu+\phi_2)+O(\mu^{2\nu_0}),\\
 z_{l-1}&=&f_{\mu}(z_1),
 \end{array}
 $$ 
where $f_{\mu}(\cdot)$ is a smooth function. With a constant $\nu_0<1$, the range of $z_1$, such that $z_{l-1}<0$, can be determined through the condition $f_{\mu}(z_1)<0$. It is represented by a union of countable disjoint positive ($\mu>0$) intervals for a fixed $\mu$ value. These intervals of one-sided $[111\cdots1]_{l}$ homoclinic orbits then correspond to new colored strips such as ones shown in Fig.~\ref{fig10}A for $B_0<1$, or Fig.~ \ref{fig11}A for $B_0>1$ that cannot overlap the strips generating  similar orbits of lower orders. The shape of the corresponding bifurcation curves in the $(\mu,\nu_0)$-parameter plane is respectively determined by the positions of the generating strips as we discussed previously. We point out that the intervals are $\mu$-value related and therefore the sides of those strips are not ``perfectly'' parallel in general. This observation, barely influencing the results, may nevertheless break some symmetry arrangements for long one-sided orbits. For example, in Figure \ref{fig10}B there may be more yellow bars on the left from the purple $\cap$-shaped bar, than on the right from it.

However, the intervals due to the condition $f_{\mu}(z_1)<0$ may merge as $\nu_0 \rightarrow 1$ provided $l>3$. As a result, a bridge can occur at a wrong position, as for example, see the bridge associated with the yellow bar (region) in Fig.~\ref{fig6}B if $l>3$. In the case $l=4$  this can be explained by solving system (\ref{ch15}) for the range of $\mu$-values for which $z_3<0$. We start off by solving the last equation to determine the range of $z_2$ such that $z_3<0$. This range is given by a union of the purple intervals on the $z_2$-axis as shown in Fig~\ref{fig14}A, where the longest interval $I$ is labeled  for further explanation. The range of $z_2$ is then used in the second last equation $z_2=\mu-B_0R^{1-\nu_0} z_1 ^{\nu_0} \sin ( \Omega_0 \ln z_1+\phi_2)+O(z_1^{2\nu_0})$ to find recursively the range of $z_1$ for which $z_3<0$. It is represented by another union of countable intervals, among which the intervals $J_1$ and $J_2$ of $z_1$ are derived from the interval $I$ of $z_2$ on one period of the sine-function of $z_1$ depicted in Panel~B of Fig.~\ref{fig14}. Finally, the range of $z_1$-values is then employed into the first equation $z_1=\mu-B_0R^{1-\nu_0}\mu^{\nu_0} \sin ( \Omega_0 \ln \mu+\phi_2)+O(\mu^{2\nu_0})$ to obtain the $\mu$-values for such $z_3<0$, which   is a union of countable $\mu$-intervals. Here,  intervals such as  $k_{11}$ and $k_{12}$ are derived from the interval $J_1$ of $z_1$ on one period of the sine-function of $\mu$ shown in  Fig.~\ref{fig14}C, while intervals such as $k_{21}$ and $k_{22}$ are derived from the interval $J_2$ of $z_1$ on the same period. In the limit $\nu_0 \rightarrow 1$, both intervals $J_1$ and $J_2$  will merge, see Fig.~\ref{fig14}B. If intervals $J_1$ and $J_2$ merge before $k_{11}$ and $k_{12}$ do, then intervals $k_{11}$ and $k_{21}$ coalesce to generate a bridge section on the bifurcation curve corresponding to [1111]-homoclinic orbit, and so do $k_{12}$ and $k_{22}$ as well, see Fig.~\ref{fig14}C.

\subsection{\bf  [100$\cdots$]--homoclinic bifurcations}
Next, consider the configuration where the 1D unstable separatrix $\Gamma_1$ first makes a single loop on one side  of the stable manifold $W^s$ , and then $l-1$ one-sided  loops ($\mu<0$) on its opposite side, before it returns to the saddle-focus. Its symbolic code is hence written as $[100\cdots\,0]_l$. 
Such an orbit may be viewed as a longer version of the inverted homoclinic connection depicted in Fig.~\ref{fig4}D.

The Poincar\'{e} map with only the dominant terms (equations for $\varphi$s are omitted) can be written as follows: 
\begin{equation}  
\begin{array}{lcl}
z_1&=&-\mu-B_0R^{1-\nu_0}(-\mu)^{\nu_0} \sin (\Omega_0 \ln (-\mu)+\phi_2)+\\&&+O((-\mu)^{2\nu_0}),\\
z_2&=&-\mu+B_0R^{1-\nu_0} (-z_1) ^{\nu_0} \sin(\Omega_0 \ln (-z_1)+\phi_2)+\\&&+O((-z_1)^{2\nu_0}),\\
 & \vdots& \\
z_{l-2}&=&-\mu(\pm B_0R^{1-\nu_0})(-z_{l-3}) ^{\nu_0} \sin(\Omega_0\ln(-z_{l-3})+\phi_2)+\\&&+O((-z_{l-3})^{2\nu_0}),\\
z_{l-1}&=&-\mu(\pm B_0R^{1-\nu_0})(-z_{l-2}) ^{\nu_0} \sin(\Omega_0\ln(-z_{l-2})+\phi_2)+\\&&+O((-z_{l-2})^{2\nu_0}), 
\end{array} 
\label{ch16}
\end{equation}
with alternating $(+/-)$ in the equations above.

To determine the bifurcation unfolding of such homoclinic orbits,  one evaluates the range, say $L_l$,  of $\mu$-values for which $z_{l-1}<0$. It is evident that $L_l \subset L_{l-1}$ because $z_{l-2}<0$ in Eqs.~(\ref{ch16}). Therefore, the set $\{L_i\}_{i=2}^{\infty}$ decreases for any $0<\nu_0<1$.  System (\ref{ch16}) can be recast as
 $$
 \begin{array}{lcl}
 z_1&= &-\mu-B_0R^{1-\nu_0}(-\mu)^{\nu_0} \sin (\Omega_0 \ln (-\mu)+\phi_2)+\\&&+O((-\mu)^{2\nu_0}),\\
  z_{l-1}&=& g_{\mu}(z_1),
  \end{array}
 $$ 
where $g_{\mu}(\cdot)$ is a smooth function. The range of $z_1$-values for which $z_{l-1}<0$ is a union of countable disjoint negative ($\mu<0$) intervals for a fixed $\mu$. As discussed previously, the bifurcation unfolding corresponding to  $[100\cdots\,0]_l$-homoclinic orbits can be illustrated using Figs.~\ref{fig12} and \ref{fig13}: the same structures re-emerge within the middle green $\cap$-shaped bar, while narrow bridges and bars re-emerge inside bridges and bars, respectively, except for   the  center green bar. However starting with $l>3$, new bridges can reside inside bars and new bars can reside within a bridge, and so forth.

\subsection{Mixed multi-loops}

Let us finally discuss mixed multi-loops  -- informally, those are longer homoclinic orbits that are neither solely left nor right sided at the end. Each such corresponding map will be a mix of equations from systems (\ref{ch15}) and (\ref{ch16}). By omitting small term $\mu$ for simplification, define these two sin-functions: 
$$\begin{array}{lcl}
f_0(x)&=&\pm B_0R^{1-\nu_0} (-x) ^{\nu_0} \sin(\Omega_0 \ln (-x)+\phi_2),\\
f_1(x)&=&-B_0R^{1-\nu_0} x ^{\nu_0} \sin(\Omega_0 \ln x+\phi_2).
\end{array}$$
Then, the map for a $[1\alpha_1 \alpha_2  \cdots  \alpha_{l-1}]$-homoclinic orbit is given by 
\begin{equation}
\begin{array}{lcl}
z_1=f_{\alpha_1}(\mu), ~~~ z_2=f_{\alpha_2}(z_1), ~~ \cdots, ~~ z_{l-2}=f_{\alpha_{l-2}}(z_{l-3}), \\
 z_{l-1}=f_{\alpha_{l-1}}(z_{l-2}), \label{ch17}
\end{array}
\end{equation}
 where $\alpha_i \in \{0,1\}$ for $1 \le i \le l$. Let us recast the above equation~(\ref{ch17})
 as
\begin{equation}
\begin{array}{lcl}
z_{l-2}=F(\mu)=f_{\alpha_{l-2}} \circ f_{\alpha_{l-3}} \circ \cdots \circ f_{\alpha_1}(\mu). \\
z_{l-1}=f_{\alpha_{l-1}}(z_{l-2}), \label{ch18}
\end{array}
\end{equation}
We showed earlier that both $\frac{df_0}{dx}, \frac{df_1}{dx}\approx \infty$, except for within small regions near their extreme points. Therefore, $\frac{dF}{d\mu} =\frac{df_{\alpha_{l-2}}}{dz_{l-3}} \cdot \frac{df_{\alpha_{l-3}}}{dz_{l-4}} \cdots \frac{df_{\alpha_1}}{d\mu} \approx \infty$ except for small regions around critical points. Thus, the graph of $F(\mu)$ is a collection of nearly vertical line segments, which is obviously bounded. It has an abundance of  extreme points because $\frac{dF}{d\mu}=0$ where $\frac{df_{\alpha_i}}{dz_{i-1}}=0$ ($1 \le i \le l-2$ and $z_0=\mu$). Figure~\ref{fig15} illustrates such oscillating graphs of the function $z_{l-2}=F(\mu)$ when $\mu>0$; it remains similar if $\mu<0$. For the sake of structural visibility, the graph per se and monotone parts are sketched not as steep as they should be in the theory. The blue projection intervals on the horizontal axis represent all $\mu$-values for which $z_{l-2}<0$; i.e., they represent parameter intervals corresponding to $[1 \alpha_1 \alpha_2  \cdots, \alpha_{l-2}]$-mixed homoclinic loops. 

If $\alpha_{l-1}=0$, then without loss of generality, we can assume that 
$$
f_{\alpha_{l-1}}(z_{l-2})=-B_0R^{1-\nu_0} (-z_{l-2}) ^{\nu_0} \sin(\Omega_0 \ln (-z_{l-2})+\phi_2).
$$ 
The solutions satisfying the condition $z_{l-1}=f_{\alpha_{l-1}}(z_{l-2})<0$ are given by $z_{l-2} \in  (-e^{-\frac{2m\pi}{\Omega_0} + \frac{\pi}{\Omega_0}-\frac{\phi_2}{\Omega_0}},-e^{-\frac{2m\pi}{\Omega_0}-\frac{\phi_2}{\Omega_0}}) = I_m$, where $m \in \mathbb{Z}$ is large enough. Let the interval between $I_m$ and $I_{m+1}$ be notated as $J_m=(-e^{-\frac{2m\pi}{\Omega_0} -\frac{2\pi}{\Omega_0}-\frac{\phi_2}{\Omega_0}},-e^{-\frac{2m\pi}{\Omega_0}+\frac{\pi}{\Omega_0}-\frac{\phi_2}{\Omega_0}})$. It is easy to see that $I_{m+1}/I_m=J_{m+1}/J_m=e^{-\frac{2\pi}{\Omega_0}}$; here,  $\{I_m\}$ are the colored strips accumulating to the $\mu$-axis from below in Fig.~\ref{fig15}, where only three such  strips are sampled. The projection intervals of the interception of the colored strips with the graph $z_{l-2}=F(\mu)$ are the sought for the $[1 \alpha_1  \alpha_2, \cdots  \alpha_{l-1}]_l$-orbit intervals on the $\mu$-axis;  they reside inside the blue intervals corresponding to the $[1 \alpha_1  \alpha_2  \cdots  \alpha_{l-2}]_{(l-1)}$-orbits. The $[1  \alpha_1 \alpha_2 \cdots \alpha_{l-2}]_{(l-1)}$-homoclinic orbits occur at the endpoints of the blue intervals, such as the points $C$ and $D$ indicated in Fig.~\ref{fig15}. Inside the blue region [CD], the green interval is much wider than the rest of the corresponding intervals because it is due to one of the extreme points of the function $z_{l-2}=F(\mu)$, whereas smaller intervals are due to nearly vertical oscillation of the given sinusoidal function.

If $\alpha_{l-1}=1$, then the solutions satisfying the inequality $z_{l-1}=f_{\alpha_{l-1}}(z_{l-2})<0$ are positive ($\mu>0$) intervals, and therefore the intervals for $[1 \alpha_1  \alpha_2  \cdots  \alpha_{l-1}]_l$-orbits are located within the gaps between the intervals for the shorter $[1 \alpha_1 \alpha_2  \cdots  \alpha_{l-2}]_{(l-1)}$-orbits, as seen from  Fig.~\ref{fig15}B. Note from this figure that the middle  interval for a $l$-long orbit is significantly wider than the rest of such ones  because it is resulted from a flatter section of the graph of the function $z_{l-2}=F(\mu)$.

It will be shown below that basic biparametric sweeps of systems with saddle-foci visibly reveal some of the largest or principle homoclinic bifurcation structures, which are likely due to such flat regions of the function $z_{l-2}=F(\mu)$, while ones due to steep oscillatory graph sections are often too narrow to be well detected and require some parameter recalling. 

In a small neighborhood of an endpoint of each interval for  $[1  \alpha_1  \alpha_2  \cdots  \alpha_{l-2}]-{(l-1)}$-orbits, such as $C$ or $D$ in Figure \ref{fig15}, the oscillating function $z_{l-2}=F(\mu)$ is nearly linear. Therefore, near the end points the parameter intervals for the $[1 \alpha_1  \alpha_2  \cdots \alpha_{l-1}]_l$-homoclinic orbits hold a scalability ratio $e^{-\frac{2\pi}{\Omega_0}}$  for both width and distance, following from the relationship:  $I_{m+1}/I_m=J_{m+1}/J_m=e^{-\frac{2\pi}{\Omega_0}}$.

The following theorem concludes our arguments and reckoning above:
\begin{theorem}
Under the conditions of Theorem~1, in the ($\mu,\,\nu_0)$-parameter space, all bifurcation curves corresponding  to $[\alpha_1 \alpha_2  \cdots \alpha_l]_l$- homoclinic orbits, where $\alpha_i \in \{0,1\}$,  $1 \le i \le l$, (i) are embedded  between the bifurcation curves corresponding to $[\alpha_1  \alpha_2  \cdots  \alpha_{(l-1)}]_{(l-1)}$-homoclinic orbits provided that $\alpha_l=0$,  or (ii) they reside within the gaps between the bifurcation curves  corresponding to $[\alpha_1 \alpha_2 \cdots  \alpha_{(l-1)} ]_{(l-1)}$-homoclinic orbits if $\alpha_l=1$. The scalability ratio for both the widths and the distances of the curves is given by $e^{-\frac{2\pi}{\Omega_0}}$, near the bifurcations of $[\alpha_1 \alpha_2  \cdots  \alpha_{(l-1)}]_{(l-1)}$-orbits.
\end{theorem}

\begin{figure*}[t!]
\begin{center}
\includegraphics[width=.99\linewidth]{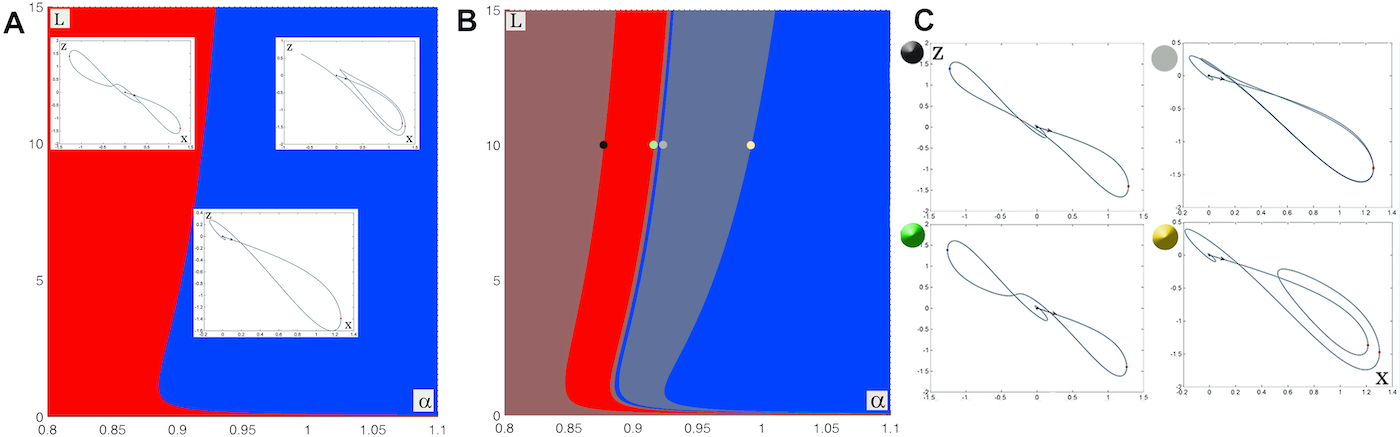}
\caption{ (A) ($\alpha, L$)-biparameteric sweep using [1-2]-long binary sequences reveals the primary homoclinic bifurcation curve [1], separating the red and blue regions (see a homoclinic trajectory in the bottom inset). The trajectories in the red region have symbolic encoding of the form $\{1,0,\ldots\}$ (see trajectory in top left inset) and those in the blue region have $\{1,1,\ldots \}$ (top right inset). (B) [1-3]-long binary sequence reveals two further homoclinic bifurcation curves, [10] (with the black and green points) and [11] (with the gray and yellow points). The positions of the black, green, gray and yellow points in the $\alpha - l$ plane are given by $(0.876898493756,9.995)$, $(0.91631772114,9.995)$, $(0.921727874,9.995)$ and $(0.991649733,9.995)$, resp., and the corresponding homoclinic orbits are shown in (C).}\label{fig17}
\end{center}
\end{figure*}

\section{Homoclinic bifurcations in the smooth Chua circuit}\label{sec:3}

Let us get back to the smooth Chua model~(\ref{ch1}) and discuss the bifurcations of its equilibria. The summary can be found in the $(a,\,b)$-bifurcation diagram shown in Fig.~\ref{fig3}. Recall that  the system~(\ref{ch1}) has three equilibrium states located at  $O(0,0,0)$, $O_1(-1,0,1)$ and $O_2(1,0,-1)$. The curve labeled by NSF, standing for ``neutral'' saddle, with the saddle index $\nu=1$ or zero saddle value $\sigma_1=0$, at the origin is given by  $b=\frac{(a^2-33a+36)(a-6)}{36(3-a)}$. The NDSF curve of $O$ is given by $a=6$;  this abbreviation stands for the saddle-focus at $0$ with zero divergency $\sigma_2=0$ or $\nu=1/2$ for saddle-foci; the sum of all three characteristic exponents is negative below this level where the system remains dissipative, whereas it becomes positive above it, making the space volume expand near the origin. While there is no curve corresponding to $\nu=0$ at the origin $O$, however, the curve for $\nu=\xi$ at $O$ is given by $b=\left (\frac{7a(a-6)}{12} -\frac{\xi(a-6)^3}{36(1-2\xi)^2} \right ) /(a\xi-3)$ for $\xi \neq 1/2$ and $0 <\xi <1$. It has an asymptote $a=3/\xi$, and therefore, we can approximate the curve $\nu=\xi$  using $a=3/\xi$. The curve $\mu=0$ defined for modeling the 1D map~(\ref{1dmap}) must be associated with the $H8$-curve in fig~\ref{fig3} that corresponds to the primary homoclinic bifurcation in the system~(\ref{ch1}), as illustrated in Fig.~\ref{fig5}. It is easy to see the correspondence between  the $(\mu,\nu_0)$-parameter plane and the $(a,\,b)$-parameter plane, so  the system (\ref{ch1}) is an ideal example to showcase the theory built in Section~2.

\subsection{Symbolic computational method}

  Bi-parametric sweepings of the system (\ref{ch1}) is done by using the computational methods originally introduced in our earlier papers~\cite{BSS12,XBS14} with a few changes, see also the following papers \cite{Barrio2013,pusuluri2017unraveling,pusuluri2018homoclinic,pusuluri2019symbolic,Pusuluri2020Chapter,pusulurihomoclinicCNSNS} 

 We follow the trajectory that initiates from the right unstable separatrix $\Gamma_1$ of the origin and record ``1'' when it loops around the right equilibrium $O_2(1,0,-1)$ and ``0'' when it loops around the left equilibrium $O_1(-1,0,1)$, see Fig.~\ref{fig2}A. Alternatively, we can use $x$-traces to convert into binary sequences so that ``1" stands for a positive maximum greater than 1, and record ``$0$'' when  $x$ reaches a negative minimum smaller than -1, see Fig.\ref{fig2}B. We can skip the very first symbol as it is always ``1'' by construction. Such a binary sequence, also knows as a kneading sequence, is recorded for a pair of $a$- and $b$-parameter values to create a bi-parametric sweep.  Next, the binary sequence is converted to a decimal number by using this following rule:
$$
K(a,b) = \sum_{n=i}^j {\kappa}_{n}\,q^{(j-n+1)} , 
$$ 
where $\{ {\kappa}_{n} \}_{n=i}^j $ is the corresponding binary sequence with ${\kappa}_{n}=\{ 0, \,1\}$, and $i,j$ are positive	 integers with $i \le j$ (the first $i-1$ binary symbols are skipped). $0<q<1$ is chosen for such a formal power series to converge. In this study, we set $q=0.5$, and keep  $K(a,b)$ in the range $[0,1]$. This decimal number is known as the kneading invariant. By construction, the $K$-values range between 0 and 1. The boundary values are set by the periodic sequences $\{\overline{0}\}$) and $\{\overline{1}\}$, respectively, for infinitely long sequences. 

\begin{figure*}[t!]
\begin{center}
\includegraphics[width=.9\linewidth]{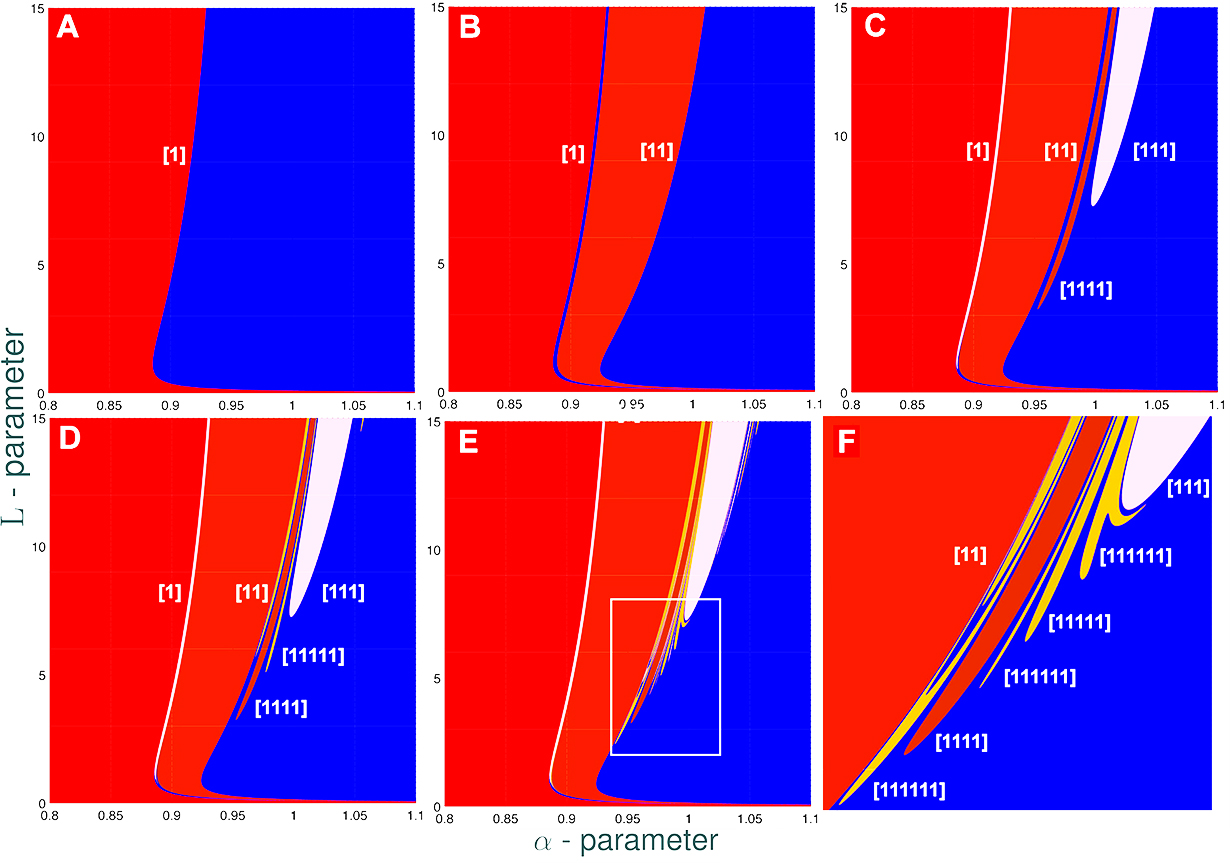}
\caption{A series of biparametric sweeps with longer binary sequences reveal an increasing hierarchy of bifurcation curves of one-sided homoclinic orbits beginning with (A) primary one; (B) two identified curves for double [11]/[00]-loops; (C) bifurcation curves for triple [111]/[000] (wide white), and [1111]/[0000]-orbits (narrow dark red)  and longer  one-sided orbits in (D)-(E). (F)  Magnified inset from (E) revealing various $\cap$-shaped and ``Y''-shaped secondary bifurcation curves (here, corresponding to [111111]-homoclinic orbits) as predicted by the theory.}\label{fig18}
\end{center}
\end{figure*}

\begin{figure*}[t!] 
	\begin{center}
		\includegraphics[width=.99\linewidth]{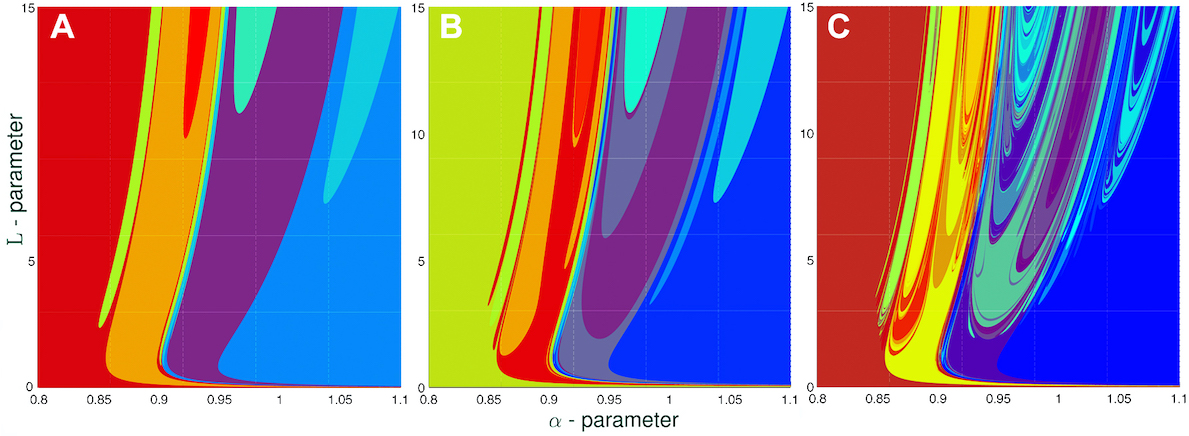}
		\caption{ Three snapshots of the $(\alpha,L)$-biparametric sweeps of the Chua system with increasing length of symbolic sequences from [2,4] in (A) through [2,5] and [2,9] in (B-D), resp., reveal the complexity and organization of mixed multi-loop homoclinic bifurcation curves of the system. 
		}\label{fig19}
	\end{center}
\end{figure*}

\begin{figure*}[t!]
	\begin{center}
		\includegraphics[width=.95\linewidth]{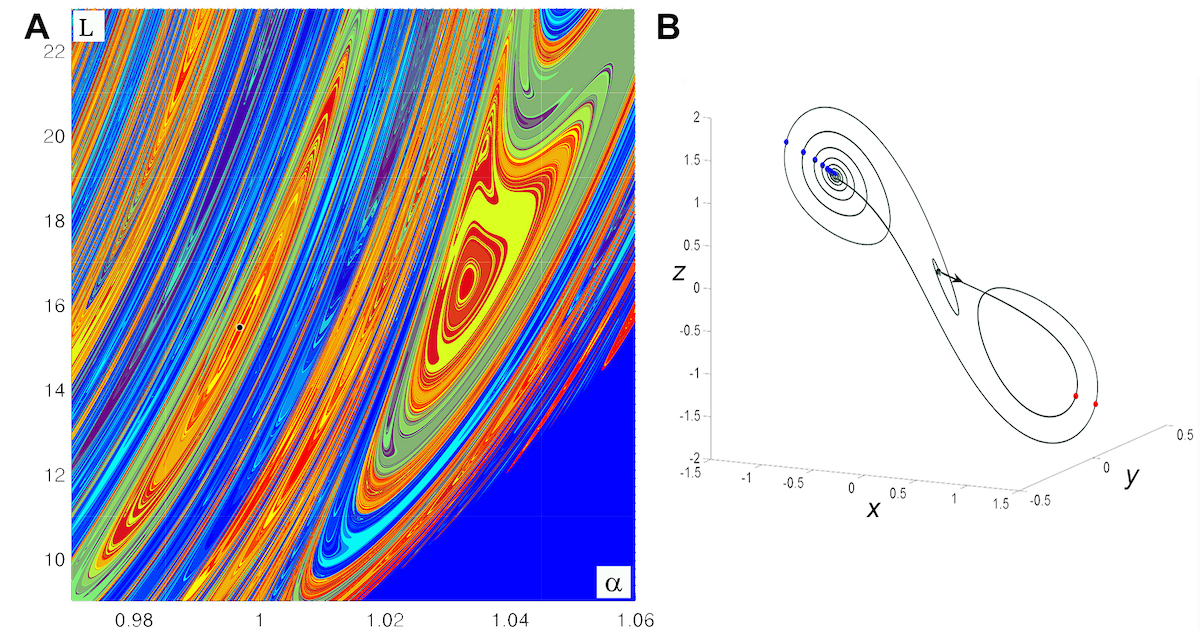}
		\caption{(A) A biparametric [7-17]-long sweep made of 42 sub-panels, each with $1000 \times 1000$ points, showing a T-point -- the black dot $(0.9971,15.2888)$ -- at the center of the characteristic spiral.  (B) A two-way heteroclinic connection at the T-point (black dot in panel~A): the separatrix $\Gamma_1$ of the saddle-focus at the origin $O$ terminates at the left saddle-focus $O_2$, whereas an outgoing trajectory spiraling away from $O_2$ converges to $O$.}\label{fig20}
	\end{center}
\end{figure*}

\begin{figure*}[t!]
	\begin{center}
		\includegraphics[width=.9\linewidth]{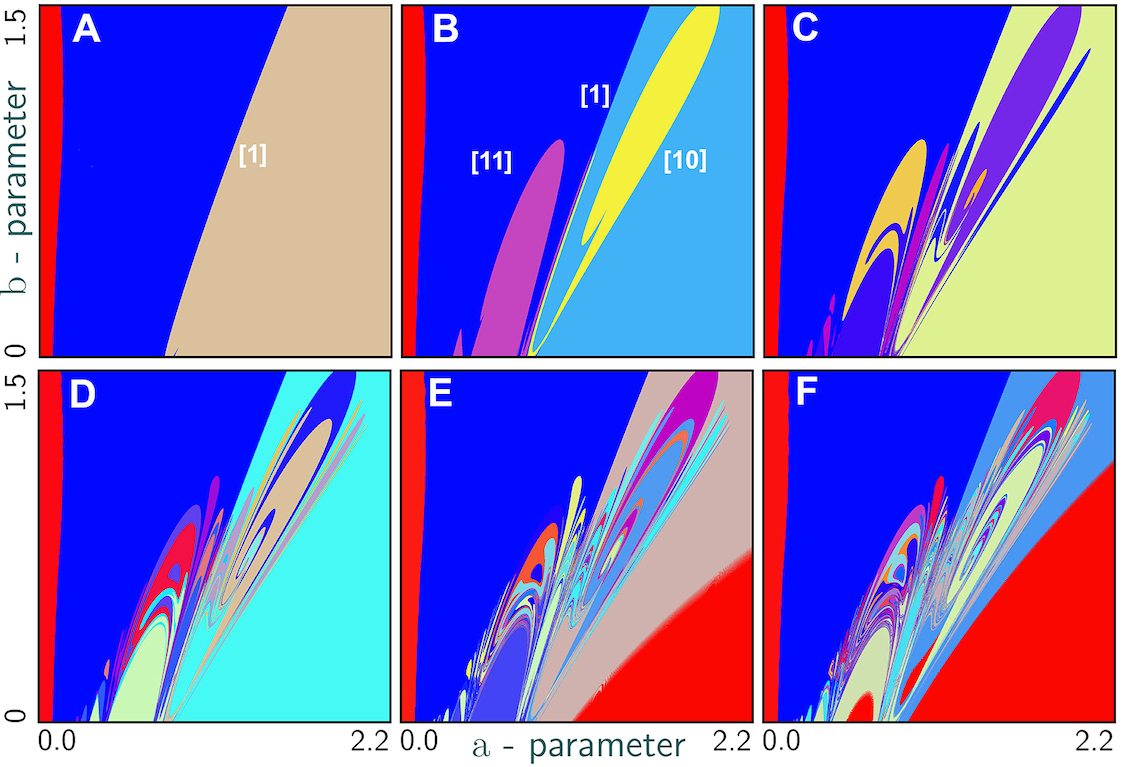}
		\caption{Bi-parametric sweeps with the increasing length of symbolic encoding, from 2 through 9, of the separatrix $\Gamma_1$ of the saddle-focus at the origin in the ACST-model using the original  $(a,b)$-parameters.  In the red regions, the 1D unstable separatrix $\Gamma_1$ of the origin escapes to infinity after a few or multiple turns around the saddle-foci $O_{1,2}$. One can observe several formed and forming T-points with characteristic spirals and nested circles around. A vicinity of the the T-point is magnified in Fig.~\ref{fig25}A below. }\label{fig21}
	\end{center}
\end{figure*}

\begin{figure*}[t!]
	\begin{center}
		\includegraphics[width=.9\linewidth]{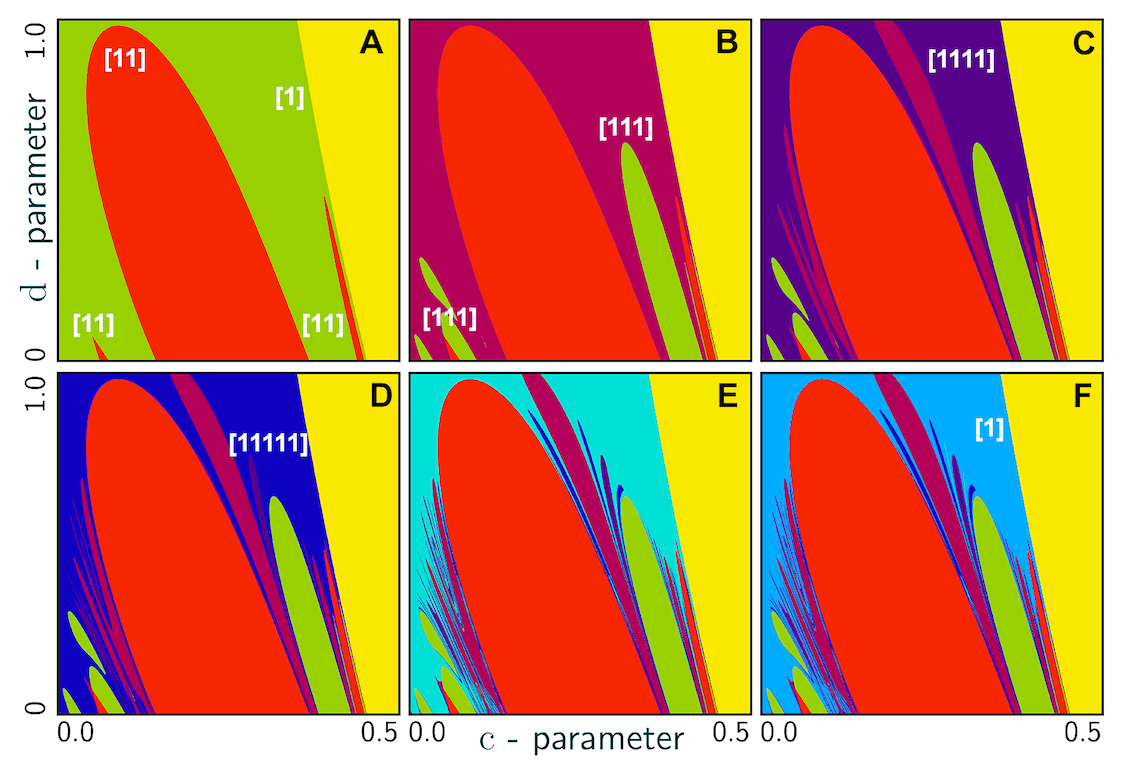}
		\caption{The biparametric sweeps of the cubic ACST-model with the new $(c,d)-$parameters  revealing the ordered intricacy of homoclinic bifurcations of one-sided homoclinic orbits encoded with [$111\cdots$], using [3,3]- through [3,8]-long symbols from (A) to (F). Panels (B-F) depict how new bifurcation curves corresponding to longer homoclinic orbits progressively fill in the gaps between ones corresponding to the shorter one-sided homoclinic orbits, to the left of the primary [1]-curve, as follows from the theory; compare with Figs.~10B-11B and the sweep for the Chua circuit in Fig.~18 above.	
			%3,3 to 3,8	
		}\label{fig22}
	\end{center}
\end{figure*}

\begin{figure*}[t!]
	\begin{center}
		\includegraphics[width=.9\linewidth]{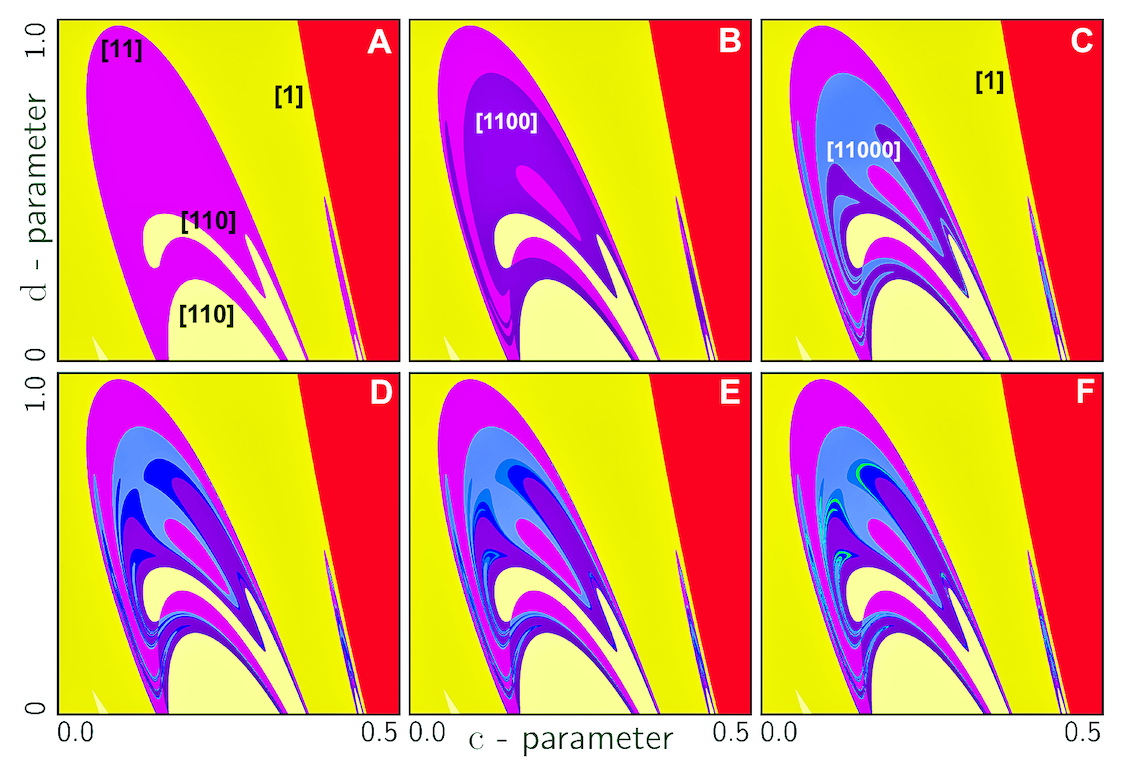}
		\caption{The biparametric sweeps of the cubic ACST-model with the new $(c,d)$-parameters  revealing how the bifurcation curves of longer ``left''-sided [$110\cdots$]  homoclinic orbits  are nested within the larger structures corresponding to double [$11$]-orbits, beginning with [110], [1100] and [11000] in (A)-(C). Panels (B-F) depicting how new left-sided bifurcation curves corresponding to longer homoclinic orbits progressively fill in the gaps between ones corresponding to the shorter one-sided homoclinic orbits below the [11]-curve, as follows from the theory; compare with the theoretical diagrams in Figs.~\ref{fig10}B-\ref{fig11}B and the bi-parametric sweep of the Chua circuit in Fig.~18 above. Here, panels (A)-(F) are obtained using symbolic sequences [3,4]- through [3,9].		%3,4 to 3,9
		}\label{fig23}
	\end{center}
\end{figure*}

\begin{figure*}[t!]
	\begin{center}
		\includegraphics[width=.9\linewidth]{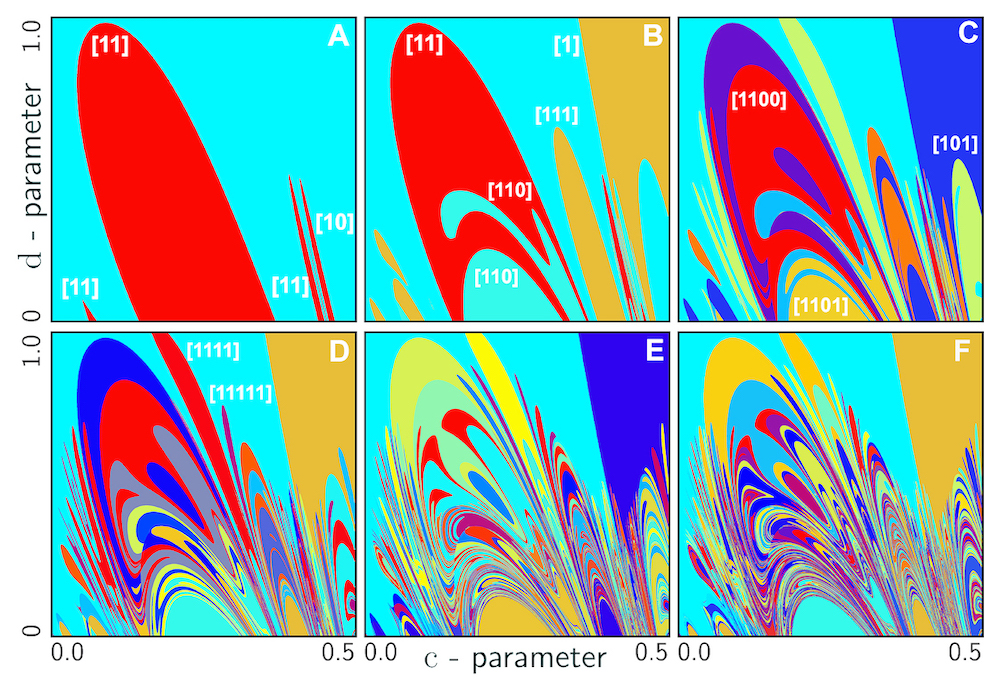}
		\caption{The $(c,d)$-parameter sweeps of the cubic ACST-model to disclose the combined bifurcation unfoldings with self-similar, scaled-down organization, due to a plethora of homoclinic orbits in its phase space; compare with Fig.~\ref{fig21} using the original parameters and the sweeps of the Chua circuit in Figs.~\ref{fig16}B and \ref{fig19}. Panels (A)-(F) are obtained using symbolic sequences [3,3]- through [3,8].
		}\label{fig24}
	\end{center}
\end{figure*}

Numerical integration is performed using a 4th-order Runge-Kutta method with a fixed step-size. The computation of trajectories across different parameter values is parallelized using GPUs. Data visualization is done in Python. A colormap takes $K$-values into $2^{8}$ discrete bins of RGB-color values, assigned from 0 through 1 for each channel of red, green and blue colors, in decreasing, random and increasing order, respectively. With such a colormap, we can assign a unique color to a single kneading invariant to produce a {\em colorful} sweep with  $1000 \times 1000$ points in the biparametric plane, as shown in Figs.~\ref{fig16}--\ref{fig26}. Parameter values that produce topologically similar trajectories result in identical sequences $\{ {\kappa}_{n} \}$, and therefore, have the same $K$-values and colors in a biparametric sweep. We employ two symbolic approaches using: (i) short, 1 $\le i \le j \le 10$, binary sequences  to detect a plethora of homoclinic bifurcations (see Figs.~\ref{fig16}--\ref{fig24}), and (ii) long, typically $600\le n \le 1000$, ones to detect stability windows within chaos-land in the parameter sweeps (see Fig.~\ref{fig26}). On a GPU-powered workstation, such a sweep takes from a few seconds to several minutes depending on the sequence length and the sweep resolution. 

Let us first discuss the first approach. By construction, a borderline between distinct colored regions in sweeps employing short sequences, is a homoclinic bifurcation curve in the parameter space. In theory, one can detect up to $2^{10}$ homoclinic bifurcations in such sweeps with $10$ binary symbols.  

Figure~\ref{fig16}A represents a short ($a,\,b$)-parameter sweep of the cubic Chua model~(\ref{ch1}). It is overlaid with the neutral saddle-focus bifurcation NSF-curve, saddle-to-saddle-focus transition (S-SF)-curve, and zero divergence NDS/NDSF-curve ($\nu=1/2$). Of our interest is the narrow wedge in the diagram that embraces many homoclinic bifurcation curves (of various colors), thus suggesting the onset of chaotic dynamics as it is located between the NSF and the S-SF curves, where the origin is the Shilnikov saddle-focus with $\nu<1$. The solid colors, blue and brown, indicate the regions of simple,  Morse-Smale dynamics with stable equilibria and/or periodic orbits, see trajectory snapshots on the  pathway $b=6$ in Fig.~\ref{fig3}.    

For a better look inside the wedge, we will apply a parameter transformation to widen this region.  First let us identify a sector in Fig. ~\ref{fig16}A, which is bounded by white curves, with its tip located  at $(1.8623, 1.8743)$. Let us introduce two new polar-coordinates/parameters: $\alpha$ and $L$; here, $\alpha$ is an angular variable inside the sector and the positive axis and $L$ is the length of the sector, i.e., radial variable. The transformation is then given by 
\begin{equation}
a= 1.8623 + L  \cos(\alpha) , \quad  b=1.8743+ L \sin(\alpha) , \label{ch19}
\end{equation}
where $\alpha \in [0.8,1.05]$ and $L \in [0,15]$. The sweep using new $(\alpha,\, L)$-parameters is shown in Figure \ref{fig16}B. It shows a plethora of bifurcation curves representing various one-sided and mixed homoclinic orbits. In what follows, we will attempt to figure out the ordered intricacy and universality of the organization of such bifurcation curves in the Chua model and other such $\mathbb{Z}_2$-symmetric systems. 

Figure~\ref{fig17}A represents the shortest sweep of length 2 to reveal the primary homoclinic bifurcation $H_1$ [1] of the saddle-focus at the origin. It occurs on the borderline of two regions: red and blue where the binary sequences start with [$10\cdots$] and [$11\cdots$], respectively. Increasing the length lets us disclose at least two pairs of bifurcation curves in the parameter diagram (Fig.~\ref{fig17}B), corresponding to the double homoclinic loops [10] and [11] (see Fig.~\ref{fig17}C for phase trajectories). Arguably, they all have a U-shape, stretched horizontally. These bifurcation curves for double loops lie on both sides of the primary one $H_1$  ($\mu=0$), as the sweep discloses. Here, the left-side of $H_1$ corresponds to $\mu<0$ and the right side corresponds to $\mu>0$, if we refer to the modeling 1D maps studied above. Figure~\ref{fig17}B detects well two principal bifurcation curves (labeled with black and yellow dots) away from $H_1$, on either sides (two other curves labeled with green and gray dots are immediately next to $H_1$, and not seen clearly). They correspond to the largest blue $\cup$-shaped bars in Figs.~\ref{fig7}B or \ref{fig10}B. The rest of the countably many curves are too close to $H_1$  to be identified in the sweep of the current scale; we discussed the reasons in Section~\uppercase\expandafter{\romannumeral2\relax} above. 

The sweep utilizing progressively longer, [2-4], binary sequences exposes bifurcation curves corresponding to various triple loops, and so on. The longer the kneading sequence used, the higher the order of homoclinic orbits and bifurcation curves that can be revealed. For example, the sweep shown in Fig.~ \ref{fig16}B utilizes symbolic subsequences  of [6-15]-range, revealing several thousands of bifurcation curves, limited due to scaling factors.

\subsection{[$111\ldots$] or one-sided homoclinic orbits and bifurcation curves}

The detailed sweep in Fig.~\ref{fig16}B discloses bifurcation curves corresponding to all homoclinic orbits of the saddle-focus, using [6-15]-long binary sequences. Next, we would like to see what the typical bifurcation unfolding of the generic Shilnikov saddle-focus may look like.  To do so, we modify the approach  to focus only on one-sided homoclinic orbits and the corresponding bifurcation curves. This is accomplished by obtaining for each ($\alpha, L$)-parameter pair, symbolic sequences containing the same symbol (say, [$111\cdots $]) until the very first occurrence of the other symbol (``0''). If the number of ``1''s in  a symbolic sequence for some  one-sided orbit is $n$, then the corresponding kneading invariant is defined as $K(\alpha, L)=n/r$, provided $r$ is the total length of such sequence.

The results of such one-sided symbolic approach are summarized in Fig.~\ref{fig18}. Using sequences up to two symbols long, we identify the primary homoclinic bifurcation shown in Fig.~\ref{fig18}A. With an additional symbol, the well-visible [11]-bifurcation curve in Figure \ref{fig18}B is revealed. This curve corresponds to the largest and the furthest blue $\cap$-shaped bar on the right in Figs.~\ref{fig7}B or \ref{fig10}B. The two white $U$- or $\cap$-shaped regions seen in Fig.~\ref{fig18}C are due to the right triple [111]-orbits that correspond to the large purple $\cap$-bar in the center of Figs. \ref{fig9}B-\ref{fig10}B, which are due to the intersections of the sine-function with the strip $S_1$ in Figs.~\ref{fig9}B-\ref{fig10}A. All other such triple orbits occur near the primary one and their bifurcation curves are hardly detected in such sweeps at the given scale, as discussed in Section~\uppercase\expandafter{\romannumeral2\relax}~C. The bifurcation curves corresponding to [1111]-homoclinic orbits fill in the gaps between those for triple orbits, see Fig.~\ref{fig18}C. This is also the case with the narrow bifurcations curves for longer orbits, see Fig.~\ref{fig18}D-E, which are harder to observe in sweeps at the given scale. To remedy this, we magnified a small region (white inset in Fig.~\ref{fig18}E) of interest near the curve terminals in Fig.~\ref{fig18}F. According to its relative position, it must be derived from the strip $S_3$ in Figure \ref{fig10}A. This magnified sweep can well depict a pair of visible yellow regions whose boundaries correspond to [111111]-homoclinic orbits. There is a single (yellow) ``Y''-shaped  region for the same [111111]-orbits near the [111] U-shaped zone, similar to ones due to the strip $S_1$  which are shown in Figs.~\ref{fig10}-\ref{fig11}. One can also see a bridge on the top of a yellow piece that was determined earlier to be derived from the strip $S_4$ in Fig.~\ref{fig10}, as was discussed in Section~\uppercase\expandafter{\romannumeral2\relax}.

\subsection{Mixed multi-loop homoclinic bifurcation curves}

Figure \ref{fig19} demonstrates four progressive $(\alpha,L)$-parameter sweeps of the Chua model with an increasing length of symbolic sequences: from [1,4] through [1,7]. Shorter sweeps are depicted in Fig.~\ref{fig17}A-B revealing the primary and double homoclinic bifurcations of the central saddle-focus in the system. One can see that the double bifurcation curves reside  on  either side of the primary one. In these parametric sweeps, the bifurcation curves for longer orbits squeeze into the gaps between or fall inside the lower order bifurcation regions. It seems likely that the most principle  bifurcation structures and boundaries, which are visible in all the sweeps, are derived from the strips that intercept the flat extreme point regions, such as the green strip in Fig.~\ref{fig10}. All other countably many curves according to the analysis done in Section~\uppercase\expandafter{\romannumeral2\relax} are too slim to be identified. Nevertheless, some interesting details can yet be pointed out. For example, the progression from Panel~A to  Panel~B Fig.~\ref{fig19} reveals two visible bridges near the bifurcation curve for double one-sided orbits, which would have to be derived from the strip $S_2$ in Figs.~\ref{fig10}--\ref{fig11}. There are also two symmetric $\cap$-shaped bars added to the picture in the left double-loop HB piece, and they must be derived from the strip $S_4$ from the same theoretical constructions. One can see many such similarities between the theoretical and computational bifurcation diagrams.

We may assume that further detailed discussion of the bifurcation unfolding, with the quickly growing complexity presented in Figs.\ref{fig19} and \ref{fig16}B, is unnecessary, or even unrealistically comprehended, as as I.~Ovsyannikov, L.P.~Shilnikov's student and co-author~\cite{OSh86e, ovsyannikov1992systems}, joked more then three decades ago: ``the saddle-focus is as inexhaustible (infinite) as the electron.''.

 \section{Bykov T-points}
  
  Figure~\ref{fig20}A presents another bi-parameter sweep of the Chua model. There are two new patterns in it that have not been well recognized in the other sweeps. The first pattern is a family of (yellow-red) nested closed circles, while the second pattern is due to several characteristic spirals stretched nearly vertically in the sweep. Such a spiral is the distinguished feature converging to the so-called T-point of codimension-two, corresponding to a two-way heteroclinic connection between a saddle and/or saddle-foci of different topological types. Such a heteroclinic connection between the central saddle-focus $O$ of (2,1)-type at the origin and the left saddle-focus $O_2$ of (1,2)-type is shown in Fig.~\ref{fig20}B. Due to the symmetry, there are always two such connections. In this phase space projection, the 1D unstable separatrix $\Gamma_1$ of $O$ and one of two 1D stable separatrices of $0_2$ coincided in the 3D phase space of the Chua model. This constitutes a one-way heteroclinic connection. Meanwhile, the 2D unstable manifold of $O_2$ and the 2D stable manifold of the origin cross along a trajectory connecting both the saddle-foci, see Fig.~\ref{fig20}. This makes this heteroclinic connection a two-way one.  According to the original research by V.~Bykov~\cite{Bykov:1999,B00}, the occurrence of a single primary $T$-point implies that there are infinitely many $T$-points nearby. A few other T-points with characteristic spirals can be also recognized in the sweep in Fig~\ref{fig20}A. Moreover, by virtue of the theory, unlike the case of the T-point of the saddle---saddle-focus connection with a single principal spiral terminating in it, the unfolding of the T-point of two saddle-foci includes two ``transverse'' spirals in the parameter plane, each representing the bifurcation curve of either Shilnikov saddle-focus with a homoclinic orbit with an incrementally increasing the number of turns around the other saddle-focus as the center of the spiral is approached with each revolution. As the Chua model is a dissipative system, we cannot employ our method to detect such bifurcation curves related to the saddle-foci $0_{1,2}$ to integrate solutions backward in time. Instead, one should use parameter continuation software, such as MatCont, that is designed to solve both initial- and boundary-value problems to continue unstable solutions. We reemphasize that according to Bykov~\cite{B80}, there are infinitely many T-points in a symmetric saddle-focus system such as the Chua model and various Lorenz-like systems, see Refs.~\cite{BSS12,XBS14,Barrio2013,pusuluri2017unraveling,pusuluri2018homoclinic}  including 3D parameter space reconstructions near T-points examined in Ref.~\cite{pusulurihomoclinicCNSNS}
 
 Concerning the nested circle pattern in Fig.~\ref{fig20}B, this happens when the bi-parameter sweep  cuts a higher, say three-dimensional parameter space not throughout a T-point (it become a space line in 3D) but only slices 2D spiraling surfaces wrapping around it at a different angle, see more in Ref.~\cite{pusulurihomoclinicCNSNS}
  This is also called a ``non-transverse'' T-point in Ref.~\cite{AM06}

\section{Symmetric ACST-mode with cubic nonlinearity}

Let us finally consider a second example to showcase the symbolic approach and to disclose the universality of homoclinic and heteroclinic patterns due to Shilnikov saddle-foci. The model (\ref{act}) is the ${\mathbb Z}_2$-symmetric extension of the generic asymptotic normal form, called the quadratic ACST-model, 
\begin{equation}\label{asct1}
\dot{x} = y, \quad \dot{y}=z, \quad \dot{z}=-y -b\,z + a\,x \left (1-x \right ), 
\end{equation}
with $(a\,b)>0$ being bifurcation parameters, describing locally  occurring  bifurcation in systems, near an equilibrium state with three zero characteristic exponents. Using the cubic term $x^3$ instead of $x^2$ in system~(\ref{asct1}) lets us worry less, computationally, about homoclinic orbits running away from the Shilnikov saddle-focus at the origin. Still, we can examine the basic homoclinic bifurcations in full generality, if we focus on one-sided orbits only in Eqs.~\ref{act}. This model also has three equilibria: the origin $O$ can be saddle-focus of the topological (2,1)-type, while $O_{1,2}(0,0, \pm 1)$ become saddle-foci of the (1,2)-type, after a supercritical Andronov-Hopf bifurcation following the period-doubling cascade initiating the onset of chaos in it, see Fig.~\ref{fig3}.     

Several snapshots of bi-parametric sweeps  with an increasing, from 2 through 8, length of symbolic encoding  of the ACST-model with the original $(a,b)$-parameters are shown in Fig.~\ref{fig21}. The area painted in red color in Panel~E is where the solutions of the model escape to infinity. One can observe from this figure that the $\cap$-shaped bifurcation curves  of longer homoclinic orbits are revealed in matching pairs in a similar fashion as the smooth Chua model. The sweeps in Fig.~\ref{fig21}D-E also disclose the location of several formed and forming T-points with their characteristic spirals and nested circles corresponding to heteroclinic connections between the saddle-foci.

Following the same approach as before, by introducing two new parameters, $c$ and $d$ with the aid of this transformation:
\begin{equation}\label{newcoord}
         a=0.24+1.76c+0.55d, \quad b=1.24c+0.81d, 
\end{equation}
we can widen up the parameter sector to provide better insights into the homoclinic unfoldings due to asymmetric one-sided or generic orbits and those due to the symmetry of the model. 

With the new parameters, the sweeps better illustrate the intrinsic organization of the Shilnikov homoclinic bifurcations in the cubic ACST-model. Let us first consider the series of sweeps shown in Fig.~\ref{fig22} representing the building hierarchy of the bifurcation unfolding representing one-sided [$11\cdots$] homoclinic orbits 
in the ${\mathbb Z}_2$-symmetric model~(\ref{act}), which would be equivalent to the generic ones occurring in system~(\ref{asct1}) with $x^2$-term. Figure~\ref{fig22}A depicts the primary [1]-homoclinic bifurcation curve and three $\cap$-shaped ones (red boundaries) corresponding to the double [11]-orbits, while Figure~\ref{fig22}B adds up several (green) curves corresponding to the triple [111]-homoclinic orbits. One can observe from the next panels in Fig.~\ref{fig22} that the new $\cap$-shaped curves for longer one-sided orbits such as [1111] and so on, fit in between the preceding ones corresponding to shorter orbits, as predicted by the theory; see Figs.~\ref{fig7}B and \ref{fig11}B above, and a similar sweep in Fig.~\ref{fig18} for the Chua circuit.

\begin{figure*}[t!]
	\begin{center}
		\includegraphics[width=.8\linewidth]{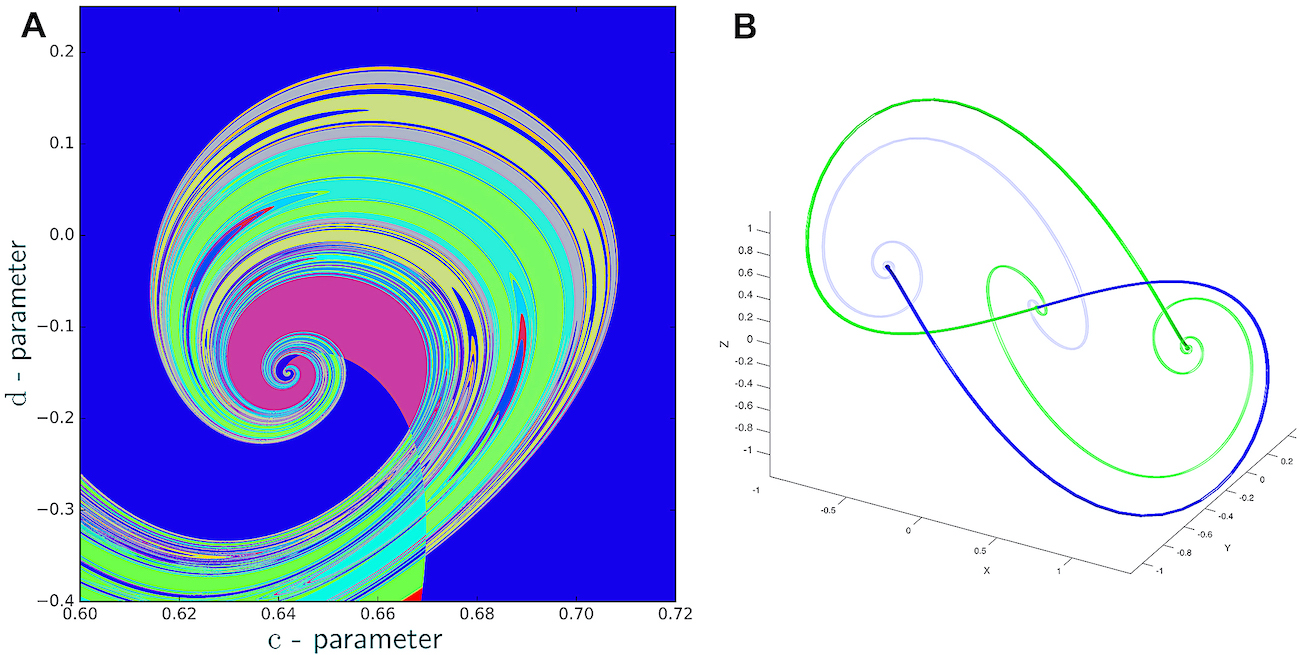}
		\caption{(A) The $(c,d)$-parameter sweep of the cubic ACST-model to reveal the primary T-point corresponding to the heteroclinic connection, [$1000\cdots \infty$] connecting the saddle-focus $O(0,0,0)$ of (2,1)-type with the saddle-foci $O_{1}$ (and $O_{2}$ due to the symmetry) of (1,2)-type in the phase space of the cubic ACST-model at $(c,d) \simeq (0.642, -0.15)$ as Panel~B depicts; compare with Fig.~\ref{fig21} using the original parameters and the sweeps of the Chua circuit in Figs.~\ref{fig16}B and \ref{fig20}. 
		}\label{fig25}
	\end{center}
\end{figure*}

Next, let us discuss how the $\cap$-shaped [11]-regions in Fig.~\ref{fig22}A are populated by bifurcation curves corresponding to left-sided homoclinic orbits encoded as [$1100\cdots$]. This is illustrated by a series of such sweeps in Fig.~\ref{fig23} of an increasing length to reveal up to six [0]s following the initial block [11]. The corresponding symbolic ranges are given by [3,4]- through [3,9].  One can see from the initial sweep in Fig.~\ref{fig23}A that the largest  [11]-region now includes two (yellow) islands whose boundaries correspond to the [110]-homoclinic orbits. The lower one is of the $\cap$-shape, while the one above is of the $Y$-shape, just like in the bifurcation sketch in Fig.~\ref{fig12}B. As the sweeping length is increased, more complex bifurcation structures for orbits such as [1100] and so on, start filling in the spaces between the matching borders for shorted orbits. We let the Reader her/himself try to figure out a self-similar order, if any, of this puzzle.  

Finally, the sweeps in Fig.~\ref{fig24} amalgamate step-by-step all identified bifurcation structures corresponding to one-sided and mixed homoclinic orbits in the cubic ACST-model. While one can easily follow   the first building steps in the Panels~\ref{fig24}A-D incorporating large structures discussed above, the last two Panels~\ref{fig24}E-F can be only inspected visually by merely stating that they incorporate smaller self-similar ones corresponding more complex orbits. In addition, one can also spot several families of nested circles due to non-transverse T-points for heteroclinic connections that interfere with our primary targets -- homoclinic bifurcations, to further recursively complicate  this overall global bifurcation unfolding, beyond feasible limits.

Figure~\ref{fig25}A, concluding this section, magnifies a vicinity of the primary Bykov T-point shown in Figs.~\ref{fig21}D-F near (0.642, -0.15) in these new ($c,d$)-parameters introduced in Eqs.~(\ref{newcoord}). Thus figure also reveals a multiplicity of secondary T-points squeezed between the spirals corresponding to homoclinic orbits with shorter encodings $[1000\cdots]$ of the origin. The demarcation curve ending at the T-point is an artifact due to the coding and color-map algorithms: with each revolution the number of zeros in the binary sequence increases incrementally by one, which makes the color change;

\section{Long-term symbolic approach to detect stability windows within chaosland}

\begin{figure*}[t!]
\begin{center}
\includegraphics[width=.99\linewidth]{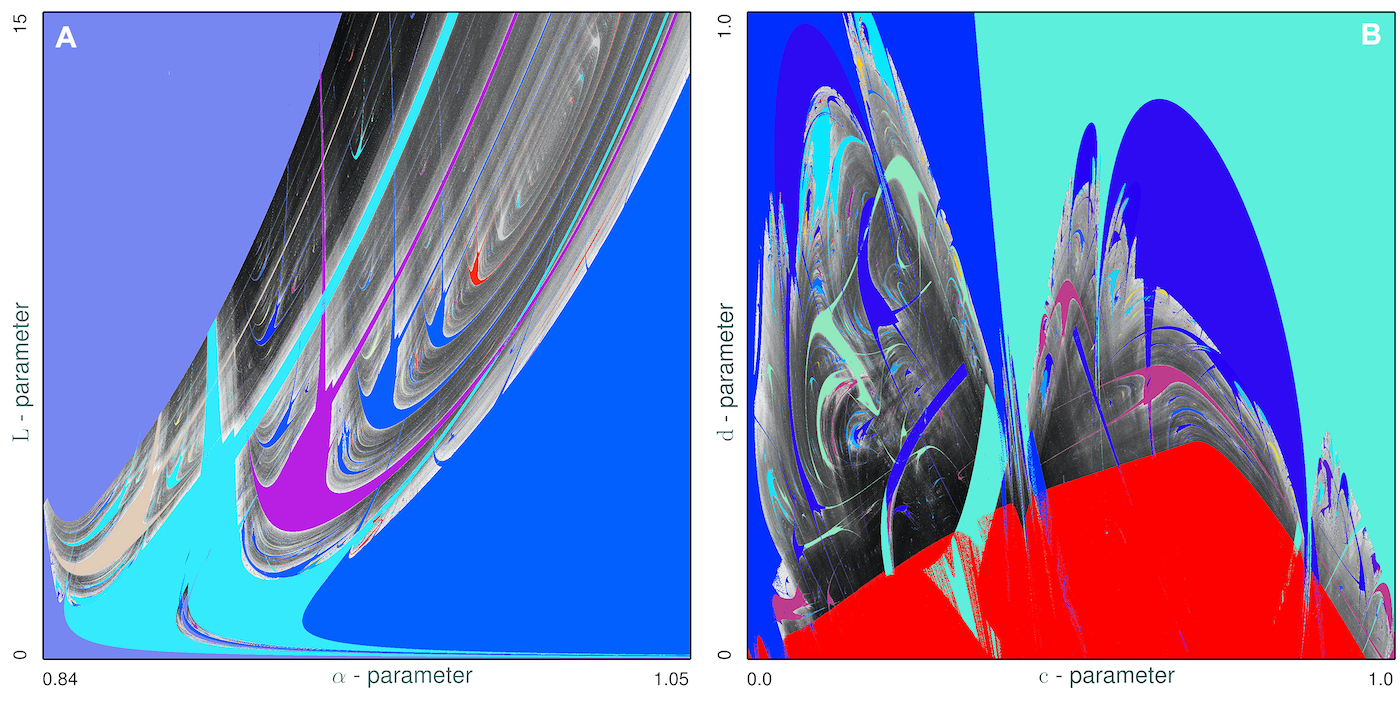}
\caption{Long [601--1000]-symbol DCP-sweeps on $5000 \times 5000$-size grid to reveal stability windows (aka ``shrimps'') (in solid colors) due to saddle-node bifurcations of periodic orbits and chaos-land (grayish regions -- with darker gray implying greater complexity) due to the Shilnikov saddle-focus in the parameter space of the Chua circuit (A) and the cubic ACST-model (B) with the transformed parameters. In the red region in (B), the solutions of the ACST-model run to infinity. }\label{fig26}
\end{center}
\end{figure*}

We have recently developed an approach called the ``Deterministic Chaos Prospector'' (DCP)  (available as an open source toolkit at  \url{https://bitbucket.org/pusuluri_krishna/deterministicchaosprospector/}) whose significance for the study of homoclinic bifurcations is that, not only can it reveal the short term transient dynamics and the underlying homoclinic, heteroclinic, saddle, and T-point spiral structures, but it can also be employed to examine the long term behavior and to detect the regions of simple dynamics due to stable equilibria and periodic orbits, and ones corresponding to chaos in the Chua, ACST and other systems such as the various Lorenz-like and R\"ossler models  \cite{pusuluri2018homoclinic,pusuluri2019symbolic,Pusuluri2020Chapter,pusulurihomoclinicCNSNS,Rossler2020}. The underlying idea is that a trajectory integrated long enough for some parameter values may  eventually converge to an exponentially stable attractor with a non-changing symbolic encoding that occupies some existence region filled out with a solid color in the parameter sweep (see Fig.~\ref{fig26}).  On the contrary, by virtue of structural instability, this is not the case for the (grayish) regions of deterministic chaotic dynamics. Such sweeps are obtained by computing trajectories from the same or different initial conditions long enough, so that after skipping some initial transients, we can still generate sufficiently long binary sequences, say of [600--1000]-length.  

The sweeps are constructed by first analyzing each long binary sequence (after omitting a transient) to detect periodicity. Periodic sequences (corresponding to simple, i.e., stable dynamics and structurally stable) of different periods are marked with different solid colors in the sweep. Aperiodic sequences (complex -- structurally-unstable dynamics) representing chaotic trajectories, are processed using the Lempel-Ziv (LZ) compression algorithm to measure their complexity \cite{lempel1976complexity}. Greater LZ-complexity indicates greater instability and is shown in darker gray. Further details of DCP can be found in Refs.\cite{pusuluri2018homoclinic, pusulurihomoclinicCNSNS, Rossler2020}

The long-term sweeps of the Chua circuit and the ACST-model are demonstrated in Figs.~\ref{fig26}A and \ref{fig26}B, respectively.  Both reveal exceptionally well, a plethora of stability windows (solid color) with  distinct periodic orbits, and regions of chaos shown in gray colors. Moreover,  we re-emphasize the darker gray pixels are associated with more developed chaos in the given models. Note that some stability windows are known as ``shrimps'' due to their shape. Such a shrimp is formed by transverse saddle-node bifurcations of periodic orbits that are typically caused by homoclinic tangencies due to spiraling saddle-foci in these and other systems. In the blue regions on opposite sides in Figs.~\ref{fig26}A, the binary sequences for the separatrix $\Gamma_1$ include periodic blocks $\{\overline 1\}$ and $\{\overline{10} \}$, respectively. While $\{\overline 1\}$ may correspond to a stable equilibrium state $O_1$, and stable periodic orbits around it or even a chaotic attractor emerging through a period doubling cascade of the former, the block $\{\overline{10} \}$ may be associated with various symmetric and asymmetric stable figure-8 periodic orbits. In Fig.~\ref{fig26}B, red color marks the region where the trajectories escape to infinity.

\section{Conclusions and discussions}

We developed a general theory of homoclinic bifurcations of the Shilnikov saddle-focus in ${\mathbb Z}_2$-symmetric systems. It discloses the ordered intricacy of corresponding structures and their organization in the bifurcation unfoldings of such systems, including  the scalability ratio $e^{-\frac{2\pi}{\Omega_0}}$ of width and distance between two sequentially close $l$- and  $(l+1)$-homoclinic orbits.   

The theoretical foundations were implemented using a novel algorithm of symbolic, binary description to examine and demonstrate the universal organization of Shilnikov homoclinic bifurcations in two symmetric systems: the smooth Chua circuit and the cubic asymptotic normal form -- the Arneodo-Coullet-Spiegel-Tresser model.    

We demonstrated how recently developed toolkit, Deterministic Chaos Prospector with GPU parallelization, can quickly reveal the regions of simple and chaotic dynamics in the parameter space of the selected models. 

The theory and the methodology created in this study can further advance new theoretical ideas and computational approaches for a better understanding of the origin and the universal structure of deterministic  chaos in full generality, including diverse application from  mathematical, physical and biological sciences.  

\section{Acknowledgement}

We are very grateful to L.P. Shilnikov for inspiration and guidance, and know that he would have warmly welcomed this paper.  
 
We thank the Brains and Behavior initiative of Georgia State University for the fellowships awarded to T. Xing and K. Pusuluri. The Shilnikov NeurDS lab thanks the NVIDIA Corporation for donating the Tesla K40 GPUs that were actively used in this study. A.~Shilnikov acknowledges a partial funding support from the Laboratory of Dynamical Systems and Applications at NRU HSE, grant No. 075-15-2019- 1931 from the Ministry of Science and Higher Education of Russian Federation.

\section*{Data/Code Availability}
The DCP code used in this study is open source and freely available at \url{https://bitbucket.org/pusuluri_krishna/deterministicchaosprospector/}.

\section*{References}
%\bibliography{saddle_focus.bib}

\providecommand{\noopsort}[1]{}\providecommand{\singleletter}[1]{#1}%

\end{document}